\documentclass[10pt]{amsart}

\newcounter{mnote}
\let\oldmarginpar\marginpar
\renewcommand\marginpar[1]{\-\oldmarginpar[\raggedleft\footnotesize #1]%
{\raggedright\footnotesize #1}}

\usepackage[english]{babel}

\usepackage{amssymb}
\usepackage{amsmath}
\usepackage{mathrsfs}
\usepackage{stmaryrd}
\usepackage{chemarrow}
\usepackage[norelsize]{algorithm2e}
\usepackage{enumerate}
\usepackage{graphicx}
\usepackage[all]{xy}
\usepackage{tikz,tikz-cd}
\usetikzlibrary{arrows}
\usepackage{multirow}
\usepackage[colorinlistoftodos]{todonotes}
\usepackage[colorlinks=true, allcolors=blue]{hyperref}
\usepackage[hyperpageref]{backref}
\usepackage{adjustbox}

\usepackage{bm}         

\usepackage{booktabs} 
\usepackage{multirow} 

\usepackage{subcaption} %
\usepackage{float}
\usepackage{subcaption}
\newtheorem{theorem}{Theorem}[section]
\newtheorem{lemma}[theorem]{Lemma}
\newtheorem{corollary}[theorem]{Corollary}
\newtheorem{proposition}[theorem]{Proposition}

\newtheorem{example}[theorem]{Example}

\newtheorem{remark}[theorem]{Remark}

\newcommand{\bs}{\boldsymbol}
\newcommand{\mcal}{\mathcal}

\newcommand{\curl}{\operatorname{curl}}
\renewcommand{\div}{\operatorname{div}}
\newcommand{\grad}{\operatorname{grad}}
\newcommand{\tr}{\operatorname{tr}}

\newcommand{\dev}{\operatorname{dev}}
\newcommand{\sym}{\operatorname{sym}}

\newcommand{\hess}{\operatorname{hess}}
\newcommand{\dg}{\operatorname{dev \, grad}}

\newcommand{\bbR}{\mathbb{R}}
\newcommand{\bbS}{\mathbb{S}}
\newcommand{\bbT}{\mathbb{T}}
\newcommand{\jump}[1]{\llbracket #1 \rrbracket}
\newcommand{\interl}[1]{\interleave #1 \interleave_{1, h} }
\newcommand{\RT}{\operatorname{RT}}
\newcommand{\Ivem}{I_h^{\mathrm{VE}}}
\newcommand{\Int}{I_h^{\mathrm{nt}}}

\newcommand{\tIvem}{I_h^{\mathrm{W}}}

\newcommand{\Qnt}{Q_h^{\mathrm{nt}}}
\newcommand{\gradvem}{\grad_{h}^{\Pi}}
\newcommand{\curlvem}{\curl_{h}^{\Pi}}

\newcommand{\innerF}{\mathring{\mathcal{F}}_h}

\newcommand{\VhRT}{V_{k,h}^{\mathrm{RT}}}
\newcommand{\Vhnt}{\mathring{\Sigma}_{k,h}^{\mathrm{nt}}}
\newcommand{\VhvemS}{\mathring{V}_{k+1,h}^{\bbS,\mathrm{VE}}}
\newcommand{\Whvem}{\mathring{W}_{k+2,h}^{\mathrm{VE}}}
\newcommand{\Whvemm}{\mathring{W}_{k+3,h}^{\mathrm{VE}}}

\newcommand{\PT}{\Pi_{k+2,T}^{H^2}}
\newcommand{\Ph}{\Pi_{k+2,h}^{H^2}}

\newcommand{\Hcd}{H^{-1}(\overline{\curl\div},\Omega;\mathbb{T})}

\begin{document}

\title[]{An Arbitrary-Order Distributional Finite Element Method for the Curl--Dev Grad Stokes System}
\author{Ziwen Gu}%
\address{School of Mathematics, Shanghai University of Finance and Economics, Shanghai 200433, China}%
\email{2024310134@stu.sufe.edu.cn}%
\author{Xuehai Huang}%
\address{School of Mathematics, Shanghai University of Finance and Economics, Shanghai 200433, China}%
\email{huang.xuehai@sufe.edu.cn}%

\thanks{The second author is the corresponding author and was supported by the National Natural Science Foundation of China (Grant No. 12671432).}


\begin{abstract}
We develop an arbitrary-order distributional finite element method for the three-dimensional curl--dev grad Stokes system governed by the Hessian complex. On tetrahedral meshes, the method couples degree-$(k+1)$ nonconforming symmetric-tensor virtual elements, degree-$k$ tangential-normal continuous traceless tensor elements, and Raviart--Thomas elements of index $k$, without additional virtual-element stabilization. Building on zero-extension curl--div theory, we derive a well-posed distributional saddle-point formulation with an $H(\operatorname{div})$ vector multiplier. A divergence-augmented discrete norm equivalence controls the full Raviart--Thomas space modulo translations and underpins mesh-uniform stability for each fixed $k\ge0$. We establish optimal $O(h^{k+1})$ energy error estimates and an optimal $O(h^{k+2})$ $L^2$ error estimate for the full symmetric-tensor approximation, under suitable regularity and load approximation assumptions. As an application, the method supplies the Stokes block of a decoupled triharmonic scheme with optimal Hessian error estimates. Numerical experiments support the predicted convergence rates for computable quantities.
\end{abstract}

\maketitle

\section{Introduction}\label{sec:introduction}
Throughout the paper, $\Omega\subset\mathbb R^3$ is a bounded contractible polyhedral domain with Lipschitz boundary. We consider the generalized tensor-valued Stokes system: for $\boldsymbol g\in L^2(\Omega;\mathbb S)$, find $\boldsymbol\sigma\in H_0^1(\Omega;\mathbb S)$, $\boldsymbol p\in H_0(\div,\Omega;\mathbb T)$, and $\boldsymbol r\in H^1(\Omega;\mathbb R^3)/\mathrm{RT}$ such that
\begin{equation}\label{generalized_stokes}
\left\{
\begin{aligned}
-\Delta\boldsymbol\sigma+\sym\curl\boldsymbol p &= \boldsymbol g &&\text{in }\Omega,\\
\curl\boldsymbol\sigma+\dev\grad\boldsymbol r &=0 &&\text{in }\Omega,\\
\div\boldsymbol p&=0 &&\text{in }\Omega.
\end{aligned}\right.
\end{equation}
We refer to \eqref{generalized_stokes} as the curl--dev grad Stokes system. Here $\mathbb S$ and $\mathbb T$ are the symmetric and traceless matrices, respectively, and $\mathrm{RT}:=\{a\boldsymbol x+\boldsymbol b:a\in\mathbb R,\ \boldsymbol b\in\mathbb R^3\}$. For $X=L^2(\Omega;\mathbb R^3)$, $H^1(\Omega;\mathbb R^3)$, or $H(\div,\Omega)$, the quotient $X/\mathrm{RT}$ is identified with its $L^2$-orthogonal representatives and carries the norm inherited from $X$.

The system \eqref{generalized_stokes} is governed by the Hessian complex \eqref{hessiancomplex} and arises in decompositions of higher-order elliptic problems \cite{ChenHuang2018,Gallistl2017}. Existing low-order discretizations include those in \cite[Section~5.3]{Gallistl2017} and the nonconforming linear method in \cite{GuHuang2026}. The latter relies on a discrete Helmholtz decomposition of piecewise-constant traceless tensors. Rather than extending that lowest-order discrete decomposition, we develop an arbitrary-order distributional method on tetrahedral meshes with a different tensor--multiplier coupling and establish its stability through a divergence-augmented norm equivalence and commuting interpolation.

Distributional formulations accommodate tensor derivatives of negative regularity through weak differential couplings. Related developments include tangential-normal stress elements for the classical Stokes equations \cite{Jay2020}, distributional Hessian and divdiv complexes \cite{HuLinZhang2025}, and distributional curl div complexes for quad-curl problems \cite{ChenHuangZhang2025SINUM}. The classical Stokes method in \cite{ChenHuangZhao2026} combines tangential-normal tensors with $H(\div)$-conforming velocities and obtains superconvergence.

Building on the zero-extension curl--div theory of \cite{Hu2025}, we establish the exact $H(\div)'$-valued realization \eqref{distributional_hessian_complex} of the Hessian complex and derive the distributional saddle-point formulation \eqref{continue_generalized_stokes}. Its vector multiplier lies in $H(\div,\Omega)/\mathrm{RT}$, rather than $L^2/\mathrm{RT}$ or $H^1/\mathrm{RT}$ as in \cite{GuHuang2026}, providing the continuous setting for a Raviart--Thomas discretization. A stable divergence lifting in $\dev\grad H^1(\Omega;\mathbb R^3)$, together with exactness, yields well-posedness and equivalence to \eqref{generalized_stokes}.

For each integer $k\ge0$, the method couples degree-$(k+1)$ nonconforming symmetric-tensor virtual elements, degree-$k$ tangential-normal continuous traceless tensors, and Raviart--Thomas elements of index $k$. Since the tangential-normal space is generally nonconforming to the continuous graph space, we discretize the divergence pairing through a weak operator. Polynomial projections make the virtual-element gradient and curl terms computable, while a projected-gradient norm equivalence avoids additional virtual-element stabilization.

The continuous vector multiplier is defined modulo $\mathrm{RT}$, whereas its discrete counterpart lies in $\VhRT/\mathbb R^3$. The latter retains the dilation field $\bs s_h(\bs x)=\bs x$, for which $\dev\grad_w\bs s_h=0$ but $\grad\bs s_h=\bs I$. The Raviart--Thomas moments support commuting interpolation, but stability is required on the full multiplier space, not only on the divergence-free subspace considered in \cite{ChenHuangZhao2026}. We prove that the divergence controls the polynomial and tangential-trace residuals associated with the Raviart--Thomas enrichment, leading to the norm equivalence
\[
\interl{\boldsymbol s_h}\simeq
\bigl(\|\dev\grad_w\boldsymbol s_h\|^2+\|\div\boldsymbol s_h\|^2\bigr)^{1/2},
\qquad \boldsymbol s_h\in\VhRT/\mathbb R^3.
\]
Here $\VhRT$ is the Raviart--Thomas space, $\dev\grad_w$ is the weak deviatoric gradient, and the broken norm includes elementwise gradients and tangential jumps. Together with a discrete inf-sup condition obtained by commuting interpolation, this equivalence establishes mesh-uniform stability for each fixed $k$. The divergence term added to the discrete bilinear form is consistent because the exact vector multiplier vanishes.

For sufficiently regular solutions and exact discrete loads, we prove optimal $O(h^{k+1})$ energy error estimates for all $k\ge0$. Under the corresponding dual regularity assumptions, the full symmetric-tensor approximation has the optimal $O(h^{k+2})$ error in $L^2$. For $k=1$, a computable piecewise-linear reconstruction of the test functions in the load functional retains this optimal rate.

As an application, we combine this Stokes method with $H^2$-nonconforming virtual elements \cite{ChenHuang2020,ChenHuang2025a} in an existing decomposition \cite{AnHuangZhang2024,ChenHuang2018,Gallistl2017,GuHuang2026} of
\begin{equation}\label{triharmonic}
-\Delta^3u=f\quad\text{in }\Omega,
\qquad u\in H_0^3(\Omega).
\end{equation}
A moment-preserving load and a polynomial tensor reconstruction preserve optimal projected and broken Hessian error estimates under the stated regularity assumptions.

The remainder of the paper is organized as follows. Section~\ref{sec:generalized_stokes} establishes the continuous distributional formulation. Section~\ref{sec:discrete_generalized_stokes} presents the discrete method, stability, and error analysis. Section~\ref{sec:triharmonic} treats the triharmonic application, Section~\ref{sec:numerical_results} reports numerical experiments, and Section~\ref{sec:conclusion} concludes the paper.

\section{A distributional formulation of the curl--dev grad Stokes system}\label{sec:generalized_stokes}
Building on the zero-extension curl--div theory of \cite{Hu2025}, this section develops a distributional saddle-point formulation of the curl--dev grad Stokes system with a vector multiplier $\bs r\in H(\div,\Omega)/\mathrm{RT}$. We construct a stable divergence lifting in $\dev\grad H^1$, establish the inf-sup stability of the resulting formulation, and prove its equivalence with \eqref{generalized_stokes}.

Let $\mathbb M$, $\mathbb S$, and $\mathbb T$ denote the spaces of $3\times3$ matrices, symmetric matrices, and traceless matrices, respectively. For $\bs\tau\in\mathbb M$, set
\[
\dev\bs\tau:=\bs\tau-\frac13\tr(\bs\tau)\bs I,
\qquad
\sym\bs\tau:=\frac12(\bs\tau+\bs\tau^{\intercal}).
\]
For a bounded Lipschitz domain $D\subset\mathbb R^3$, we use the standard Sobolev spaces $H^m(D)$ with norm $\|\cdot\|_{m,D}$ and seminorm $|\cdot|_{m,D}$, abbreviating $\|\cdot\|_{0,D}$ by $\|\cdot\|_D$. Vector- and tensor-valued spaces are written as $B(D;\mathbb X):=B(D)\otimes\mathbb X$. The space $H_0^m(D;\mathbb X)$ is the closure of $C_0^\infty(D;\mathbb X)$ in $H^m(D;\mathbb X)$. When $D=\Omega$, domain subscripts on norms, seminorms, and $L^2$ inner products are omitted. We write $\bs n$ for the outward unit normal.

The operators $\grad$ (or $\nabla$), $\curl$, and $\div$ act row-wise on tensor fields. We use
\begin{align}
\label{def:spaces}
H^{-1}(\curl\div,D;\mathbb T)
&:=\{\bs\tau\in L^2(D;\mathbb T):\div\bs\tau\in H_0(\div,D)'\},\\
\notag
H(\div,D;\mathbb T)
&:=\{\bs\tau\in L^2(D;\mathbb T):\div\bs\tau\in L^2(D;\mathbb R^3)\},\\
\notag
H_0(\div,D;\mathbb T)
&:=\{\bs\tau\in H(\div,D;\mathbb T):\bs\tau\bs n=0\text{ on }\partial D\}.
\end{align}
Here $H'$ denotes the dual of $H$. 
The notation \(A\lesssim B\) means \(A\le CB\) with a generic constant \(C\), and \(A\simeq B\) means \(A\lesssim B\) and \(B\lesssim A\).

\subsection{Zero-extension curl--div graph space and its \texorpdfstring{$H(\div)'$}{H(div)} characterization}\label{sec:zero_extension}
On the contractible domain $\Omega$, the following Hessian complex \cite{GuHuang2026} is exact:
\begin{equation}\label{hessiancomplex}
0\xrightarrow{\subset} H_0^3(\Omega)\xrightarrow{\nabla^2} H_0^1(\Omega;\mathbb S)\xrightarrow{\curl} H_0(\div,\Omega;\mathbb T)\xrightarrow{\div}L^2(\Omega; \mathbb{R}^3) / \textnormal{RT} \xrightarrow{}0.
\end{equation}
To relax the regularity requirement $ \bs  p\in H_0(\div,\Omega;\mathbb T)$ , we allow its divergence to take values in $H(\div,\Omega)'$ rather than $L^2(\Omega;\mathbb R^3)$. This leads to a vector multiplier in $H(\div,\Omega)/\mathrm{RT}$. However, the conventional curl–div space in \eqref{def:spaces} only guarantees $\div \bs p\in H_0(\div,\Omega)'$, which does not provide the required pairing. We therefore retain the boundary contribution through zero extension and construct the corresponding tensor space.

The standard distributional divergence of a function in $L^2(\Omega;\mathbb{R}^3)$ is defined by
\begin{equation*}
\langle \div \bs u,  v \rangle := -( \bs u, \nabla v), \quad \forall\,v \in C_0^{\infty}(\Omega),
\end{equation*}
and $ \div \bs u \in H^{-1} := H_0^1(\Omega)' $. To obtain an element of $H^1(\Omega)'$, we define the zero-extension divergence $ \overline{\div}$ by
\begin{equation}
\langle \overline{\div} \bs u,  v \rangle := -( \bs u, \nabla v), \quad \forall\,v \in C^{\infty}(\bar{\Omega}). \nonumber
\end{equation}
By the Cauchy--Schwarz inequality and the density of $C^{\infty}(\bar\Omega)$ in $H^1(\Omega)$, this definition extends uniquely to an element of $H^1(\Omega)'$. For $\bs u\in H(\div,\Omega)$, the Green formula gives
\[
\langle\overline{\div}\bs u,v\rangle
=(\div\bs u,v)-\langle\bs u\cdot\bs n,v\rangle_{\partial\Omega},
\qquad v\in H^1(\Omega),
\]
where the boundary term is the $H^{-1/2}(\partial\Omega)$--$H^{1/2}(\partial\Omega)$ duality pairing. Equivalently, $\overline{\div} \bs u$ is the distributional divergence of the zero extension of $\bs u$ to $\mathbb R^3$. Retaining this normal-flux contribution gives the characterization
\[
\{\bs u\in L^2(\Omega;\mathbb R^3):\overline{\div}\bs u\in L^2(\Omega)\}
=H_0(\div,\Omega),
\]
see \cite[Theorem~6.2]{Kurula2012duality}.
For a tensor $\bs p\in L^2(\Omega;\mathbb T)$, set
\begin{align*}
\langle\overline{\div}\bs p,\bs v\rangle
&:=-(\bs p,\grad\bs v),
&&\bs v\in C^\infty(\overline\Omega;\mathbb R^3),\\
\langle\overline{\curl\div}\bs p,\bs v\rangle
&:=\langle\overline{\div}\bs p,\curl\bs v\rangle
=-(\bs p,\grad\curl\bs v),
&&\bs v\in C^\infty(\overline\Omega;\mathbb R^3).
\end{align*}
We define the zero-extension curl--div graph space by
\begin{equation*}
\Hcd
:=\left\{\bs p\in L^2(\Omega;\mathbb T):
\overline{\curl\div}\bs p\in H^1(\Omega;\mathbb R^3)'\right\},
\end{equation*}
with graph norm
\begin{equation}\label{norm:Hcd_native}
\|\bs p\|_{-1,\overline{\curl\div}}^2
:=\|\bs p\|^2
+\|\overline{\curl\div}\bs p\|_{(H^1)' }^2.
\end{equation}

The following characterization follows from \cite[Lemma~5.3]{Hu2025}.
\begin{theorem}\label{thm:Hcds_Hcdss}
For $\bs p\in L^2(\Omega;\mathbb T)$,
\begin{equation}\label{char:Hcd_div}
\bs p\in\Hcd
\quad\Longleftrightarrow\quad
\overline{\div}\bs p\in H(\div,\Omega)'.
\end{equation}
Consequently,
\begin{equation}\label{eq:Hcd_norm_equivalence}
\|\bs p\|_{-1,\overline{\curl\div}}
\simeq
\left(\|\bs p\|^2+\|\overline{\div}\bs p\|_{H(\div)'}^2\right)^{1/2}.
\end{equation}
\end{theorem}
\begin{proof}
By \cite[Lemma~5.3]{Hu2025}, a functional $F\in H^1(\Omega;\mathbb R^3)'$ extends continuously to $H(\div,\Omega)$ if and only if its zero-extension curl belongs to $H^1(\Omega;\mathbb R^3)'$, with
\[
\|F\|_{H(\div)'}\simeq
\bigl(\|F\|_{(H^1)'}^2+\|\overline{\curl} F\|_{(H^1)'}^2\bigr)^{1/2}.
\]
Here $\langle\overline{\curl} F,\bs v\rangle:=\langle F,\curl\bs v\rangle$ for smooth $\bs v$ on $\overline\Omega$. Applying this characterization to $F=\overline{\div}\bs p$, and using $\|\overline{\div}\bs p\|_{(H^1)'}\lesssim\|\bs p\|$, proves \eqref{char:Hcd_div} and \eqref{eq:Hcd_norm_equivalence}.
\end{proof}

\subsection{An \texorpdfstring{$H(\div)'$}{H(div)}-realization of the Hessian complex and a stable \texorpdfstring{$\dev\grad$}{dev-grad} lifting}

\begin{theorem}\label{thm:distributional_hessian_exact}
Denote the annihilator of $\RT$ in $H(\div,\Omega)'$ as 
\begin{equation*}
Q := \{ F^* \in H(\div,\Omega)': \, \langle F^*, \bs r \rangle \, = 0, \,  \forall\,\bs r \in \RT \} \subset H(\div,\Omega)'.
\end{equation*}
Then the following \(H(\div)'\)-realization of the Hessian complex is exact:
\begin{equation}\label{distributional_hessian_complex}
0\longrightarrow H_0^3(\Omega) \xrightarrow{\ \hess \ }
H_0^1(\Omega;\mathbb{S}) \xrightarrow{\ \curl\ }
\Hcd \xrightarrow{\ \overline{\div}\ }
Q
\longrightarrow 0.
\end{equation}
Moreover, $\overline{\div}:\Hcd\to Q$ admits a stable right inverse: for every $F^*\in Q$, there exists $\bs p\in\dev\grad H^1(\Omega;\mathbb R^3) \cap \Hcd$ such that
\begin{equation}\label{stable_right_inverse}
\overline{\div}\bs p=F^*,
\qquad
\|\bs p\|_{-1,\overline{\curl\div}}
\lesssim \|F^*\|_{H(\div)'}.
\end{equation}
\end{theorem}
\begin{proof}
For $\bs p\in\Hcd$ and $\bs r\in\RT$,
$\langle\overline{\div}\bs p,\bs r\rangle=-(\bs p,\dev\grad\bs r)=0$, so $\overline{\div}\Hcd\subset Q$. Moreover,
$\langle\overline{\div}\curl\bs\tau,\bs v\rangle=-(\curl\bs\tau,\grad\bs v)=0$
for $\bs\tau\in H_0^1(\Omega;\mathbb S)$ and smooth $\bs v$, while $\curl\hess v=0$. Hence \eqref{distributional_hessian_complex} is a complex.

Exactness at $H_0^3(\Omega)$ and $H_0^1(\Omega;\mathbb S)$ follows directly from \eqref{hessiancomplex}. By the exactness of the Hessian complex~\eqref{hessiancomplex}, $\operatorname{curl}H_0^1(\Omega;\mathbb S) = H_0(\operatorname{div},\Omega;\mathbb T)\cap\ker(\operatorname{div}).$ Thus, to prove exactness at \(\Hcd\), it suffices to establish
\begin{equation}\label{div_equivalence}
\Hcd\cap\ker(\overline{\operatorname{div}})
=
H_0(\operatorname{div},\Omega;\mathbb T)\cap\ker(\operatorname{div}).
\end{equation}
If \(\bs p\in\Hcd\) satisfies \(\overline{\operatorname{div}}\bs p=0\), the characterization in Section~\ref{sec:zero_extension}, applied row-wise, implies that \(\bs p\in H_0(\operatorname{div},\Omega;\mathbb T)\) with \(\operatorname{div}\bs p=0\). Conversely, if \(\bs p\in H_0(\div,\Omega;\mathbb T)\) and \(\div\bs p=0\), Green's formula gives \(\overline{\div}\bs p=0\). Consequently, \(\overline{\curl \div} \bs  p=0\), and \(\boldsymbol p\in\Hcd\). This proves the equality and hence exactness at \(\Hcd\).

It remains to show that every \(F^*\in\mathcal Q\) can be written as \(\overline{\operatorname{div}}\boldsymbol p\), with \(\boldsymbol p\) satisfying \eqref{stable_right_inverse}. The deviatoric-gradient Korn inequality \cite[Lemma~3.2(iii)--(iv)]{PaulyZulehner2020} 
\[
\|\bs v\|_1\lesssim\|\dev\grad\bs v\|,
\qquad \bs v\in H^1(\Omega;\mathbb R^3)/\RT,
\]
implies that $\dev\grad H^1(\Omega;\mathbb R^3)$ is closed in $L^2(\Omega;\mathbb T)$. For $F^*\in Q$, define a linear functional \(\ell\) on \(\dev \grad H^1(\Omega;\mathbb R^3)\) by
\[
\ell(\dev\grad\bs v):=-\langle F^*,\bs v\rangle.
\]
Since $\ker(\dev\grad)=\RT$ and $F^*$ annihilates $\RT$, the functional is well defined, and
\[
|\ell(\dev\grad\bs v)|
\lesssim\|F^*\|_{H(\div)'}\|\dev\grad\bs v\|.
\]
By the Riesz representation theorem, there exists $\bs p\in\dev\grad H^1(\Omega;\mathbb R^3)$ such that
\begin{equation}\label{tmp_eq}
\langle F^*,\bs v\rangle=-(\bs p,\dev\grad\bs v),
\qquad \bs v\in H^1(\Omega;\mathbb R^3),
\end{equation}
and $\|\bs p\|\lesssim\|F^*\|_{H(\div)'}$. Since $\bs p$ is traceless, \eqref{tmp_eq} shows that $\overline{\div}\bs p=F^*$ on $H^1(\Omega;\mathbb R^3)$. By the density of \(H^1(\Omega;\mathbb R^3)\) in \(H(\div,\Omega)\),  $F^*$ is the unique continuous extension of $\overline{\div}\bs p$ to \(H(\operatorname{div},\Omega)\); by Theorem~\ref{thm:Hcds_Hcdss} and the preceding Riesz estimate, \(\boldsymbol p\in\Hcd\) and satisfies~\eqref{stable_right_inverse}. 
The Riesz construction is linear in \(F^*\) and therefore defines a bounded right inverse.
\end{proof}

\begin{remark}\label{rem:relation_Hu2025}
In the notation of \cite[Section~1.1]{Hu2025},
\[
\widetilde H^{-1}(\Omega)
:=\overline{C_c^\infty(\Omega)}^{\,H^{-1}(\mathbb R^3)}
\simeq H^1(\Omega)',
\]
where the test functions in the closure are extended by zero to $\mathbb R^3$. This differs from $H^{-1}(\Omega)=H_0^1(\Omega)'$. If $\widetilde{\bs p}$ is the zero extension of $\bs p$, the zero-extension curl--div space defined in \cite[Section~4]{Hu2025} is characterized by $\curl\div\widetilde{\bs p}\in\widetilde H^{-1}(\Omega;\mathbb R^3)$, and hence coincides algebraically and topologically with $\Hcd$ equipped with \eqref{norm:Hcd_native}.

Generic continuous right inverses are available in \cite[Lemma~3.16 and Theorem~3.21]{Hu2025}. The tensor \(\boldsymbol p\) constructed in Theorem~\ref{thm:distributional_hessian_exact} is the unique element of \(\operatorname{dev}\operatorname{grad}H^1(\Omega;\mathbb R^3)\) satisfying \(\overline{\operatorname{div}}\boldsymbol p=F^*\). Moreover, it has the smallest \(L^2\) norm among all tensors in \(\Hcd\) satisfying the same equation. Indeed, if $\bs z\in\Hcd$ satisfies $\overline{\div}\bs z=F^*=\overline{\div}\bs p$ and $\bs p=\dev\grad\bs v$, then
\[
(\bs z-\bs p,\bs p)
=-\langle\overline{\div}(\bs z-\bs p),\bs v\rangle=0,
\qquad
\|\bs z\|^2=\|\bs p\|^2+\|\bs z-\bs p\|^2.
\]

\end{remark}

\subsection{Distributional formulation for the curl--dev grad Stokes problem}
The distributional weak formulation of \eqref{generalized_stokes} is: find
\[
(\boldsymbol\sigma,\boldsymbol p,\boldsymbol r)
\in H_0^1(\Omega;\mathbb S)\times\Hcd\times H(\div,\Omega)/\mathrm{RT}
\]
such that
\begin{subequations}\label{continue_generalized_stokes}
\begin{align}
a(\boldsymbol{\sigma} , \boldsymbol{r}; \boldsymbol{\tau} , \boldsymbol{s}) + b( \boldsymbol{\tau}, \boldsymbol{s}; \boldsymbol{p}) & = \left( \boldsymbol{g}, \boldsymbol{\tau} \right), && \forall\, \boldsymbol{\tau} \in H_0^1(\Omega; \mathbb{S}), \boldsymbol{s} \in H(\div,\Omega) / \textrm{RT}, \label{continue_generalized_stokes1} \\
b( \boldsymbol{\sigma}, \boldsymbol{r}; \boldsymbol{q}) & = 0, && \forall\,\boldsymbol{q} \in \Hcd, \label{continue_generalized_stokes2}
\end{align}
\end{subequations}
where the bilinear forms are defined by
\begin{equation*}
a(\boldsymbol{\sigma}, \boldsymbol{r}; \boldsymbol{\tau}, \boldsymbol{s}) = ( \grad \boldsymbol{\sigma}, \grad \boldsymbol{\tau}), \quad 
b(\boldsymbol{\tau}, \boldsymbol{s}; \boldsymbol{q}) = ( \curl \boldsymbol{\tau}, \boldsymbol{q}) - \langle \overline{\div} \bs q, \bs s \rangle.
\end{equation*}
The duality pairing in $b$ is well defined on $H(\div,\Omega)/\mathrm{RT}$, since $\overline{\div}\bs q\in Q$ for every $\bs q\in\Hcd$.
The Cauchy--Schwarz inequality and \eqref{eq:Hcd_norm_equivalence} give the continuity bounds
\begin{align*}
|a(\bs\sigma,\bs r;\bs\tau,\bs s)|
&\lesssim \|\bs\sigma\|_1\|\bs\tau\|_1,\\
|b(\bs\tau,\bs s;\bs q)|
&\lesssim \bigl(\|\bs\tau\|_1+\|\bs s\|_{H(\div)}\bigr)
\|\bs q\|_{-1,\overline{\curl\div}}.
\end{align*}

\begin{lemma}\label{lem:coer_continue}
For $ \bs \tau \in H_0^1(\Omega;\mathbb{S}), \bs s \in H(\div,\Omega)/\RT$ satisfying $b(\bs \tau, \bs s; \bs q) = 0 $ for all $ \bs q \in \Hcd $, we have 
\begin{equation}\label{coer_continue}
\| \bs \tau \|^2_1 + \| \bs s \|^2_{H(\div)} \lesssim a(\bs \tau, \bs s ; \bs \tau, \bs s ).  
\end{equation}
\end{lemma}
\begin{proof}
Taking $\bs q=\curl\bs\tau$ gives $\curl\bs\tau=0$. Define
\[
F_{\boldsymbol s}(\boldsymbol v):=(\boldsymbol s,\boldsymbol v),
\qquad \boldsymbol v\in H(\operatorname{div},\Omega),
\]
then this functional is continuous on \(H(\operatorname{div},\Omega)\). Since \(\boldsymbol s\perp_{L^2}\mathrm{RT}\), we further have \(F_{\boldsymbol s}\in Q\). By Theorem~\ref{thm:distributional_hessian_exact} there exists \(\boldsymbol q_{\boldsymbol s}\in\Hcd\) such that $ \overline{\operatorname{div}}\boldsymbol q_{\boldsymbol s}=F_{\boldsymbol s}$, and it holds that
\[
0=b(\boldsymbol\tau,\boldsymbol s;\boldsymbol q_{\boldsymbol s})
=-\langle F_{\boldsymbol s},\boldsymbol s\rangle
=-\|\boldsymbol s\|^2.
\]
Therefore \(\boldsymbol s=0\), and the Poincaré inequality yields \eqref{coer_continue}. 
\end{proof}

\begin{lemma}\label{lem:inf_sup_continue}
The following inf-sup condition holds:
\begin{equation}\label{inf_sup_continue}
\| \bs p\|_{-1,\overline{\curl\div}} \lesssim
\sup_{\substack{
\bs \tau \in H_0^1(\Omega;\mathbb{S}) \\
\bs s \in  H(\div,\Omega)/\RT
}} \dfrac{b(\bs \tau,\bs s;\bs p)}{ \| \bs \tau \|_1 + \| \bs s \|_{H(\div)} }, \quad \forall \bs p \in \Hcd.
\end{equation}
\end{lemma}
\begin{proof}
For $\bs p\in\Hcd$, let $S(\bs p)$ denote the supremum on the right-hand side of \eqref{inf_sup_continue}, taken over nonzero pairs. Since $\overline{\div}\bs p\in Q$, the quotient-dual norm relation and the choice $\bs\tau=0$ give
\[
\|\overline{\div}\bs p\|_{H(\div)'}\lesssim S(\bs p).
\]
Theorem~\ref{thm:distributional_hessian_exact} provides a lifting $\bs q\in\Hcd$ in $\dev\grad H^1(\Omega;\mathbb R^3)$ such that
\[
\overline{\div}\bs q=\overline{\div}\bs p,
\qquad \|\bs q\|_{-1,\overline{\curl\div}}\lesssim S(\bs p).
\]
Set $\bs z=\bs p-\bs q$. By the exactness of \eqref{distributional_hessian_complex} and the stable curl lifting in \cite[Lemma~2]{GuHuang2026}, there exists $\bs\eta\in H_0^1(\Omega;\mathbb S)$ satisfying
\[
\curl\bs\eta=\bs z,
\qquad \|\bs\eta\|_1\lesssim\|\bs z\|.
\]
By Remark~\ref{rem:relation_Hu2025}, $(\bs z,\bs q)=0$, and hence
\[
\|\bs z\|^2=b(\bs\eta,0;\bs p)
\le S(\bs p)\|\bs\eta\|_1
\lesssim S(\bs p)\|\bs z\|.
\]
Consequently, $\|\bs z\|\lesssim S(\bs p)$. Since $\overline{\div}\bs z=0$, its graph norm equals its $L^2$ norm, and therefore
\[
\|\bs p\|_{-1,\overline{\curl\div}}
\le\|\bs z\|+\|\bs q\|_{-1,\overline{\curl\div}}
\lesssim S(\bs p).
\]
This proves \eqref{inf_sup_continue}.
\end{proof}
\begin{theorem}\label{thm:gStokes_continue_wellposed}
The weak formulation \eqref{continue_generalized_stokes} is well-posed and is equivalent, at the level of solutions, to the curl--dev grad Stokes system \eqref{generalized_stokes} in the distributional sense. Moreover, $\boldsymbol{r}=0$, $\curl\boldsymbol{\sigma}=0$, $\div\boldsymbol{p}=0$, and $\boldsymbol{p}\in H_0(\div,\Omega;\mathbb T)$, with
\begin{equation}\label{eq:continuous_stability}
\|\bs\sigma\|_1+\|\bs p\|
\lesssim\|\bs g\|_{H^{-1}(\Omega;\mathbb S)}
\lesssim\|\bs g\|.
\end{equation}
\end{theorem}
\begin{proof}
By the continuity bounds, the coercivity estimate \eqref{coer_continue}, and the inf-sup condition \eqref{inf_sup_continue}, the Babu\v{s}ka--Brezzi theory \cite{BoffiBrezziFortin2013} implies well-posedness and stability in the product norm for \eqref{continue_generalized_stokes}, which immediately yields \eqref{eq:continuous_stability}.

Applying the kernel argument of Lemma~\ref{lem:coer_continue} to \eqref{continue_generalized_stokes2} gives $\curl\boldsymbol{\sigma}=0$ and $\boldsymbol r=0$; in particular, $\bs r\in H^1(\Omega;\mathbb R^3)/\RT$, as required in \eqref{generalized_stokes}. Choosing $\boldsymbol\tau=0$ in \eqref{continue_generalized_stokes1} and testing with every $\boldsymbol s\in H(\div,\Omega)/\RT$ gives $\overline{\div}\boldsymbol p=0$. Together with \eqref{div_equivalence}, this yields $\boldsymbol p\in H_0(\div,\Omega;\mathbb T)\cap\ker(\div)$. Taking $\boldsymbol s=0$ in \eqref{continue_generalized_stokes1}, we obtain
\[
(\grad\boldsymbol\sigma,\grad\boldsymbol\tau)
+(\curl\boldsymbol\tau,\boldsymbol p)
=(\boldsymbol g,\boldsymbol\tau),
\qquad \boldsymbol\tau\in H_0^1(\Omega;\mathbb S).
\]
By the distributional definition $\langle\sym\curl\bs p,\bs\tau\rangle=(\bs p,\curl\bs\tau)$ for $\bs\tau\in H_0^1(\Omega;\mathbb S)$, this identity is equivalent to
\[
-\Delta\boldsymbol\sigma+\sym\curl\boldsymbol p=\boldsymbol g
\quad\text{in }H^{-1}(\Omega;\mathbb S).
\]
Together with $\curl\boldsymbol\sigma+\dev\grad\boldsymbol r=0$ and $\div\boldsymbol p=0$, this recovers \eqref{generalized_stokes}.

Conversely, for a solution of \eqref{generalized_stokes} in the stated spaces, $\bs p\in H_0(\div,\Omega;\mathbb T)$ and $\div\bs p=0$ imply $\overline{\div}\bs p=0$ and $\bs p\in\Hcd$. Since $\bs r\in H^1(\Omega;\mathbb R^3)/\RT$, for every $\bs q\in\Hcd$,
\[
-\langle\overline{\div}\bs q,\bs r\rangle
=(\bs q,\dev\grad\bs r),
\qquad
b(\bs\sigma,\bs r;\bs q)
=(\curl\bs\sigma+\dev\grad\bs r,\bs q)=0,
\]
which gives \eqref{continue_generalized_stokes2}; the first distributional equation and $\overline{\div}\bs p=0$ give \eqref{continue_generalized_stokes1}, proving the equivalence.
\end{proof}

Although the exact multiplier satisfies \(\boldsymbol r=0\), it serves to enforce the constraint \(\overline{\operatorname{div}}\boldsymbol p=0\) in the saddle-point formulation.

\section{Distributional finite element method for the curl--dev grad Stokes system}\label{sec:discrete_generalized_stokes}
We now develop an arbitrary-order distributional finite element method based on \eqref{continue_generalized_stokes} and establish its stability and error estimates.

\subsection{Discrete spaces}
Let $\{\mathcal T_h\}_{h>0}$ be a shape-regular family of tetrahedral meshes of $\Omega$, with $h:=\max_{T\in\mathcal T_h}h_T$. 
For each fixed polynomial order $k$, all hidden constants below are independent of the mesh size \(h\). Let \(\mathcal F_h\), \(\mathcal E_h\), and \(\mathcal V_h\) denote the sets of mesh faces, edges, and vertices; $\mathring{\mathcal{F}}_h, \mathring{\mathcal{E}}_h$ and $ \mathring{\mathcal{V}}_h $ denote their interior counterparts; and write $\mathcal F(T)$, $\mathcal E(T)$, and $\mathcal V(T)$ for those belonging to $T$. For a domain or mesh entity $D$, let $\mathbb P_k(D)$ be the polynomials of total degree at most $k$, with $\mathbb P_k(D)=\{0\}$ for $k<0$. For each face \(F\), let \(\boldsymbol n_F\) be a unit normal vector, chosen outward from \(\Omega\) on boundary faces. Let \(\boldsymbol n_{\partial T}\) denote the outward unit normal on \(\partial T\). We will abbreviate $\bs n_{\partial T}$ as $\bs n$ when the element is clear from context. In local formulas, $\partial_n$ denotes differentiation in the direction $\bs n_{\partial T}$. For an interior face $F=T_1\cap T_2$, set
\[
\jump{w}:=(w|_{T_1})|_F\bs n_F\cdot\bs n_{\partial T_1}
+(w|_{T_2})|_F\bs n_F\cdot\bs n_{\partial T_2},
\]
and set $\jump{w}:=w|_F$ on boundary faces. 
For each face \(F\), let
\[
\mathscr{T}^F:=\{\bs v\in\mathbb R^3:\bs v\cdot\bs n_F=0\},
\qquad
\Pi_F:=I-\bs n_F\bs n_F^{\intercal}
\]
denote the tangential vector space and the orthogonal tangential projection, respectively.
We use
\begin{align*}
H^1(\mathcal T_h)&:=\{v\in L^2(\Omega):v|_T\in H^1(T), \ \forall\,T\in\mathcal T_h\},
\end{align*}
and denote the elementwise differential operators by $\grad_h\,(\nabla_h)$, $\curl_h$, and $\div_h$. For $v\in H^1(\mathcal T_h)$, define the broken $H^1$ seminorm by
\[
|v|_{1,h}^2:=\sum_{T\in\mathcal T_h}|v|_{1,T}^2=\|\nabla_h v\|^2,
\]
with the same componentwise convention for vector and tensor fields.
For ease of presentation, we will use the notation
\begin{equation*}
(u, v)_{\mathcal{T}_h}:=\sum_{T \in \mathcal{T}_h} (u, v)_{T},\quad (u, v)_{\partial \mathcal{T}_h}:=\sum_{T \in \mathcal{T}_h} (u, v)_{\partial T},
\end{equation*}
\begin{equation*}
(u, v)_{\mathcal{F}_h}:=\sum_{F \in \mathcal{F}_h} (u, v)_{F},\quad (u, v)_{\innerF}:=\sum_{F \in \innerF} (u, v)_{F},
\end{equation*}
and \(\boldsymbol{n}\) in a \(\partial\mathcal T_h\)-pairing is understood elementwise as the outward normal \(\boldsymbol{n}_{\partial T}\).

Define broken polynomial spaces $\mathbb P_{k}(\mathcal{T}_h), \mathbb P_{k}(\mathcal{F}_h)$ as
\begin{equation} 
\begin{aligned}
\mathbb P_{k}(\mathcal{T}_h)  & = \{v \in L^2(\Omega) :v|_T\in \mathbb{P}_k(T), \forall\,T \in \mathcal{T}_h \} , \\
\mathbb P_{k}(\mathcal{F}_h)  & = \{v \in L^2(\mathcal{F}_h) :v|_F\in \mathbb{P}_k(F), \forall\,F \in \mathcal{F}_h \}. \nonumber
\end{aligned}
\end{equation}
Let $Q_{k,T}$ and $Q_{k,F}$ denote the $L^2$ projections onto $\mathbb P_k(T)$ and $\mathbb P_k(F)$, respectively, with the same notation used componentwise for vector- and tensor-valued fields. We also define $ Q_{k,T} =0$ for $k<0$. Let $Q_{k,\mcal T_h}$ and $Q_{k,\mcal F_h}$ denote their elementwise and facewise counterparts. Throughout this paper, $Q_{k,h}:=Q_{k,\mcal T_h}$.

For an integer $k\ge0$, we use the degree-$(k+1)$ nonconforming virtual element of \cite{ChenHuang2020,AyusoLipnikovManzini2016}, with local space
\begin{equation*}
V_{k+1}^{\mathrm{VE}}(T)
:=\{v\in H^1(T):\ \Delta v\in \mathbb P_{k-1}(T),\ \partial_n v|_F\in \mathbb P_k(F)\ \text{for all }F\in\mathcal F(T)\}.
\end{equation*}
For $k=0$, $V_1^{\mathrm{VE}}(T)=\mathbb P_1(T)$, recovering the Crouzeix--Raviart (CR) element with the face moments \eqref{dof_vem1}; see \cite[Section~3.2]{HuangTang2025}. The degrees of freedom (DoFs) are
\begin{subequations}\label{dof_vem}
\begin{align}
(v,q)_F,\quad &q\in \mathbb P_k(F),\; F\in\mathcal F(T),\label{dof_vem1}\\
(v,q)_T,\quad &q\in \mathbb P_{k-1}(T).\label{dof_vem2}
\end{align}
\end{subequations}
Define the global nonconforming virtual element (VE) space by
\begin{equation}\label{global_vem_space}
\mathring{V}_{k+1,h}^{\mathrm{VE}}
:=\{v\in L^2(\Omega): v|_T\in V_{k+1}^{\mathrm{VE}}(T) \textrm{ for all } T\in\mathcal T_h, \textrm{ and } Q_{k,\mcal F_h}\jump{v}=0\}.
\end{equation}
We set $ \VhvemS := \mathring{V}_{k+1,h}^{\mathrm{VE}} \otimes\bbS$ for symmetric tensors. The shifted indexing is intentional: the degree-$(k+1)$ virtual space provides computable projected gradient and curl operators of degree $k$, compatible with the degree-$k$ tangential-normal tensor space and the Raviart--Thomas space of index $k$ used below.

For $\bs\tau\in\VhvemS$, set
\begin{equation*}
\bs\tau_T:=Q_{k-1,h}\bs\tau,\qquad
\bs\tau_F:=Q_{k,\mcal F_h}\bs\tau,
\end{equation*}
and we use the same subscript convention for other tensor fields whenever these projections are defined. These projected element and face data are computable  and uniquely determined by the DoFs \eqref{dof_vem} of the virtual tensor $\bs\tau$.

We use the tangential-normal continuous tensor element in \cite{Jay2020,ChenHuangZhao2026}, with local shape-function space $\mathbb P_k(T;\mathbb T)$ and DoFs
\begin{subequations}\label{dof_nt}
\begin{align}
( \Pi_F \bs p \bs n_F, \bs q )_F, \quad & \bs q \in \mathbb{P}_k(F; \mathscr{T}^F),  \label{dof_nt1} \\
( \boldsymbol{p}, \boldsymbol{q} )_T, \quad & \boldsymbol{q} \in \mathbb{P}_{k-1}(T;\mathbb{T}), \label{dof_nt2}
\end{align}
\end{subequations}
and the global finite element space is 
\begin{equation*} 
\Vhnt = \{ \boldsymbol{p} \in \mathbb{P}_k(\mathcal{T}_h;\mathbb{T}): \jump{\Pi_F\bs p \bs n_F} =0 \text{ on each } F\in\mathcal{F}_h\}.
\end{equation*}
We also use the Raviart--Thomas (RT) element \cite{RaviartThomas1977,Nedelec1980,ChenHuang2022b}, with local shape-function space
\begin{equation}
V_k^{\mathrm{RT}}(T) = \mathbb{P}_k(T; \mathbb{R}^3) + \mathbb{H}_k(T) \boldsymbol{x},  \nonumber
\end{equation}
where $\mathbb{H}_k(T)$ denotes the homogeneous polynomials of degree $k$, the DoFs are
\begin{subequations}\label{dof_hdiv}
\begin{align}
\label{dof_hdiv1}
(\boldsymbol{v}\cdot\bs n_F, q)_F, \quad & q \in \mathbb{P}_k(F),\; F\in\mathcal F(T), \\
\label{dof_hdiv2}
(\boldsymbol{v}, \boldsymbol{q})_T, \quad & \boldsymbol{q} \in \mathbb{P}_{k-1}(T;\mathbb{R}^3),
\end{align}
\end{subequations}
and the global finite element space is 
\begin{equation*}
\VhRT = \{ \boldsymbol{v} \in H(\div,\Omega): \boldsymbol{v}|_T \in V_k^{\mathrm{RT}}(T) \text{ for all } T \in \mathcal{T}_h\}.
\end{equation*}
We adopt $\VhvemS$, $\Vhnt$, and $\VhRT$ to discretize $\bs\sigma$, $\bs p$, and $\bs r$, respectively.

\begin{remark}\label{rem:nt_nonconforming}
For \(\boldsymbol p_h\in\Vhnt\), tangential-normal continuity and the homogeneous boundary condition cancel the tangential terms arising from elementwise integration by parts. The remaining face terms pair the normal-normal jumps on interior faces and the normal-normal traces on boundary faces with \(\boldsymbol s\cdot\boldsymbol n_F\), where \(\boldsymbol s\in H(\operatorname{div},\Omega)\). The normal trace of \(\boldsymbol s\) on each element boundary belongs to \(H^{-1/2}(\partial T)\). For these face terms to define a bounded functional on \(H(\operatorname{div},\Omega)\), the normal-normal jumps and boundary traces must therefore satisfy the \(H^{1/2}\) compatibility conditions relating adjacent faces along their common edges, see \cite[Chapter~I]{GiraultRaviart1986}. Tangential-normal continuity and facewise polynomiality do not guarantee these conditions, cf. \cite[Remark~3]{PechsteinSchoeberl2011}. Thus \(\Vhnt\not\subset\Hcd\) in general, and we replace the continuous divergence pairing by a discrete weak divergence in the next subsection.
\end{remark}

\subsection{Weak operators and norm equivalence on Raviart--Thomas spaces}
Motivated by the preceding Remark~\ref{rem:nt_nonconforming}, we define the discrete weak divergence operator $ \div_w : \Vhnt \to \VhRT $ as
\begin{equation*}
(\div_w\bs p_h,\bs s)
:=(\div_h\bs p_h,\bs s)-
(\bs n^{\textup T}\bs p_h\bs n,
 \bs s\cdot\bs n)_{\partial \mathcal{T}_h},
\qquad \bs s\in\VhRT.
\end{equation*}
Similarly, define the weak deviatoric gradient $\dev\grad_w\colon\VhRT\to\Vhnt$ as the negative $L^2$-adjoint of $\div_w$,
\begin{equation}\label{eq_weak_devgrad}
\begin{aligned}
( \dev\grad_w \boldsymbol{r}_h, \boldsymbol{q} ): &= (\boldsymbol{q}, \dev\grad_h \boldsymbol{r}_h) - ( \Pi_F \boldsymbol q\bs n_F, \llbracket \Pi_F \boldsymbol{r}_h \rrbracket )_{\mathcal{F}_h} \\
&= -( \div_h \boldsymbol{q},  \boldsymbol{r}_h) + ( \jump{\bs n_F^{ \textup{T}} \bs q \bs n_F}, \bs r_h \cdot \bs n_F)_{\mathcal{F}_h},
\quad \bs q \in \Vhnt. 
\end{aligned}
\end{equation}
The second representation follows from elementwise integration by parts, using the normal continuity of $\bs r_h\in\VhRT$ and the tangential-normal continuity and homogeneous tangential-normal boundary condition of $\bs q\in\Vhnt$.

For $\bs\sigma_h\in\VhvemS$, define the computable $L^2$-projected gradient and curl operators by
\begin{equation*}
\gradvem\bs\sigma_h := Q_{k,h}\grad_h\bs\sigma_h, \qquad \curlvem\bs\sigma_h := Q_{k,h}\curl_h\bs\sigma_h,
\end{equation*}
which take values in \(\mathbb P_k(\mathcal T_h;\mathbb S \otimes \mathbb R^3)\) and \(\mathbb P_k(\mathcal T_h;\mathbb T)\), respectively, and can be computed on each element from the degrees of freedom \eqref{dof_vem} using the following identities:
\begin{align*}
( \gradvem \boldsymbol{\sigma}_h, \boldsymbol{\Psi} )_T  &= -(Q_{k-1,T}\boldsymbol{\sigma}_h, \div \boldsymbol{\Psi})_T + (Q_{k,F}\boldsymbol{\sigma}_h, \boldsymbol\Psi\bs n)_{\partial T}, \;\; \boldsymbol{\Psi} \in \mathbb{P}_{k}(T;\mathbb{S} \otimes \mathbb{R}^3),  \\
( \curlvem \boldsymbol{\sigma}_h,  \boldsymbol{q} )_T &= (Q_{k-1,T}\bs\sigma_{h}, \sym \curl_h \boldsymbol{q} )_T + (Q_{k,F}\bs\sigma_{h}, \boldsymbol{q} \times \bs n)_{\partial T},  \;\; \boldsymbol{q} \in \mathbb{P}_{k}(T;\mathbb{T}).
\end{align*}
By \cite[Lemma 3]{HuangTang2025}, the projected gradient satisfies the following norm equivalence on \(\VhvemS\):
\begin{equation}\label{norm_sigma}
| \bs \sigma_h |_{1,h}^2 \eqsim \| \gradvem \bs \sigma_h \|^2 , \quad \bs \sigma_h \in  \VhvemS.
\end{equation}
The broken Poincar\'e inequality \cite{Brenner2003} and the defining weak continuity in \eqref{global_vem_space} show that $ | \cdot |_{1,h} $ defines a norm on $\VhvemS$. 
Define a broken seminorm on $H^1(\mathcal T_h;\mathbb R^3)$
\begin{equation}\label{norm_s}
\interl{\bs s}^2 := \|\grad_h\bs s\|^2+\sum_{F\in\mathring{\mathcal F}_h}h_F^{-1}\|\jump{\Pi_F\bs s}\|_F^2, \quad \bs s\in H^1(\mathcal T_h;\mathbb R^3). 
\end{equation}
We identify \(\VhRT/\mathbb R^3\) with the subspace of \(\VhRT\) consisting of vector fields \(L^2\)-orthogonal to constant vectors. On \(\VhRT\), the seminorm \(\interl{\cdot}\) vanishes only on constant vectors and therefore defines a norm on this subspace. For \(k\ge1\), the trace inequality, projection approximation estimates, and the elementwise norm equivalence \(\|\nabla\boldsymbol\tau\|_T\simeq\|Q_{k,T}\nabla\boldsymbol\tau\|_T\) give the following estimate \cite[Lemma~4]{HuangTang2025}:
\begin{equation}\label{eq:vem_trace_estimate}
\|Q_{k-1,T}\bs\tau-Q_{k,F}\bs\tau\|_{F}
\lesssim h_T^{1/2}\|Q_{k,T}\grad\bs\tau\|_T,
\quad \forall\,\bs\tau\in V_{k+1}^{\mathrm{VE}}(T)\otimes\bbS, F\in\mathcal F(T).
\end{equation}

The next lemma shows that the divergence controls the deviatoric-gradient and tangential-trace residuals introduced by the Raviart--Thomas enrichment.
\begin{lemma}
For every integer $k\ge0$, $T\in\mathcal T_h$, and $\bs s\in V_k^{\mathrm{RT}}(T)$, the following estimate holds:
\begin{equation}\label{RT_estimate}
\|(I-Q_{k-1,T})\dev\grad\bs s\|_T + \sum_{F \in \mathcal{F}(T)} h_F^{-1/2}\| (I - Q_{k,F} ) \Pi_F \bs s \|_F \lesssim \|\div\bs s \|_T.
\end{equation}
\end{lemma}

\begin{proof}
Let $\bs x_T$ be the barycenter of $T$, set $\bs y=\bs x-\bs x_T$, and define
\[
\bs w:=\frac{\bs y}{k+3}(I-Q_{k-1,T})\div\bs s.
\]
Writing $\bs s=\bs p+\bs y\phi$, with $\bs p\in\mathbb P_k(T;\mathbb R^3)$ and $\phi$ homogeneous of degree $k$ in $\bs y$, Euler's identity gives $\div\bs s=\div\bs p+(k+3)\phi$. Hence
\[
\bs s-\bs w=\bs p+\bs yQ_{k-1,T}\phi\in\mathbb P_k(T;\mathbb R^3),
\qquad
\|\bs w\|_T\le\frac{h_T}{k+3}\|\div\bs s\|_T.
\]
Both projection residuals in \eqref{RT_estimate} vanish on $\mathbb P_k(T;\mathbb R^3)$. Since $\bs w\in\mathbb P_{k+1}(T;\mathbb R^3)$, the polynomial inverse and trace inequalities therefore yield
\[
\begin{aligned}
&\|(I-Q_{k-1,T})\dev\grad\bs s\|_T
+\sum_{F\in\mathcal F(T)}h_F^{-1/2}\|(I-Q_{k,F})\Pi_F\bs s\|_F\\
= ~&\|(I-Q_{k-1,T})\dev\grad\bs w\|_T
+\sum_{F\in\mathcal F(T)}h_F^{-1/2}\|(I-Q_{k,F})\Pi_F\bs w\|_F\\
\lesssim ~&h_T^{-1}\|\bs w\|_T\lesssim\|\div\bs s\|_T.
\end{aligned}
\]
The constant depends only on $k$ and the shape regularity.
\end{proof}

The unisolvence of the DoFs \eqref{dof_nt}, together with a scaling argument, yields
\begin{equation}\label{norm_equivalence_nt}
\|\bs p\|_T^2
\simeq
\|Q_{k-1,T}\bs p\|_T^2
+\sum_{F\subset\partial T} h_F\|\Pi_F\bs p \bs n_F\|_F^2 \quad\forall\,\bs p\in\mathbb P_k(T;\mathbb T).
\end{equation}

\begin{theorem}\label{lem:RT_equivalent_norm}
For every integer $k\ge0$, we have
\begin{equation}\label{eq_rt_equivalence}
\interl{\bs s_h}\simeq\|\bs s_h\|_{R,h},\qquad \bs s_h\in\VhRT,
\end{equation}
where $\|\bs s_h\|_{R,h}^2:=\|\dev\grad_w\bs s_h\|^2+\|\div\bs s_h\|^2$.
These seminorms define equivalent norms on $\VhRT/\mathbb R^3$, and
\[
\ker(\dev\grad_w)\cap\ker(\div)=\mathbb R^3
\quad\text{on }\VhRT.
\]
\end{theorem}
\begin{proof}
Choose $\bs p_h\in\Vhnt$ by prescribing its interior-face and cell DoFs as
\begin{align*}
( \Pi_F \bs p_h \bs n_F, \bs q )_F &= -h_F^{-1}( Q_{k,F}\jump{ \Pi_F \bs s_h}, \bs q )_F, \quad && \bs q \in \mathbb P_k(F; \mathscr{T}^F),\; F\in \innerF ,   \\
( \boldsymbol{p}_h, \boldsymbol{q} )_T &= ( Q_{k-1,T}\dev\grad\bs s_h, \boldsymbol{q} )_T, \quad && \boldsymbol{q} \in \mathbb P_{k-1}(T;\mathbb{T}).
\end{align*}
Summing the local norm equivalence \eqref{norm_equivalence_nt} over the mesh gives
\begin{equation}\label{eq:20261005}
\|\bs p_h\|^2 \eqsim E^2 + J^2,
\end{equation}
where
\begin{equation*}
E^2:= \sum_{ T \in \mcal T_h}
\|Q_{k-1, T} \dg \bs s_h \|_T^2, \quad
J^2: =\sum_{ F \in \innerF } h_F^{-1}\| Q_{k,F} \jump{\Pi_F \bs s_h }\|_F^2.
\end{equation*}
For $k=0$, the cell moments are absent and $E=0$. By the definition of $\dev\grad_w$ and by the choice of $\bs p_h$,
\begin{align*}
(\dev\grad_w\bs s_h, \bs p_h)  =  \,& ( \bs p_h, \dev\grad_h \bs s_h ) - ( \Pi_F \bs p_h \bs n_F, \llbracket \Pi_F \bs s_h \rrbracket )_{\innerF}   \\
= \,& E^2 + J^2 
+((I-Q_{k-1,h})\dev\grad_h\bs s_h, \bs p_h)_{\mcal T_h}. 
\end{align*}
This together with the Cauchy--Schwarz inequality, \eqref{RT_estimate} and \eqref{eq:20261005} gives
\begin{equation*}
E + J\lesssim  \|\dev\grad_w\bs s_h \|+\|\div\bs s_h \|.
\end{equation*}
Decomposing the deviatoric gradient and the face jumps into their projected and residual parts, using \(\grad_h\bs s_h=\dev\grad_h\bs s_h+\frac13(\div\bs s_h)I\), and applying \eqref{RT_estimate} elementwise and on both sides of each interior face, we obtain $\interl{\bs s_h}\lesssim\|\bs s_h\|_{R,h}$.

For the reverse inequality, let $ \bs p_h = \dg_w \bs s_h $. By the definition of $\dev\grad_w$, the Cauchy--Schwarz inequality, the discrete polynomial trace inequality, and the local norm equivalence \eqref{norm_equivalence_nt}, we have
\begin{align*}
\| \dg_w \bs s_h \|^2 & = ( \dg_w \bs s_h, \bs p_h ) \\
&= ( \bs p_h, \dev\grad_h \bs s_h ) - ( \Pi_F \bs p_h \bs n_F, \llbracket \Pi_F \bs s_h \rrbracket )_{\innerF}  \\
&\lesssim \| \bs p_h \| \interl{ \bs s_h}
= \| \dg_w \bs s_h \| \interl{ \bs s_h},
\end{align*}
which combined with $ \| \div \bs s_h \| \lesssim \| \grad_h \bs s_h \| $ gives
\(
\| \bs s_h \|_{R,h} \lesssim \interl{\bs s_h}.
\)
\end{proof}

\subsection{Discrete formulation}
Based on the weak formulation \eqref{continue_generalized_stokes}, we propose the distributional finite element method for the curl--dev grad Stokes system: Find $( \boldsymbol{\sigma_h}, \boldsymbol{p_h}, \boldsymbol{r_h}  )  \in \VhvemS \times \Vhnt \times \VhRT/\mathbb{R}^3   $ such that 
\begin{subequations}\label{discrete_generalized_stokes}
\begin{align}
a_h(\boldsymbol{\sigma}_h , \boldsymbol{r}_h; \boldsymbol{\tau} , \boldsymbol{s}) + b_h( \boldsymbol{\tau}, \boldsymbol{s}; \boldsymbol{p}_h) & = \left( \boldsymbol{g_h}, \tilde{\bs \tau} \right), && \forall\, \boldsymbol{\tau} \in  \VhvemS, \boldsymbol{s} \in \VhRT/\mathbb{R}^3 , \label{discrete_generalized_stokes1} \\
b_h( \boldsymbol{\sigma}_h, \boldsymbol{r}_h; \boldsymbol{q}) & = 0, && \forall\,\boldsymbol{q} \in \Vhnt. \label{discrete_generalized_stokes2}
\end{align}
\end{subequations}
Here $\bs g_h\in L^2(\Omega;\mathbb S)$ is an approximation of $\bs g$, and the reconstruction of test function $\widetilde{\bs\tau}$ is defined in~\eqref{eq:load_reconstruction}. The bilinear forms are
\begin{align*}
a_h( \bs \sigma_h , \bs r_h; \boldsymbol{\tau} , \boldsymbol{s}) &:= ( \gradvem \bs \sigma_h, \gradvem \boldsymbol{\tau} ) + (\div \bs r_h , \div \bs s), \\
b_h(\boldsymbol{\tau}, \boldsymbol{s}; \boldsymbol{p}_h) &:= ( \curlvem \boldsymbol{\tau} + \dev\grad_w \boldsymbol{s}, \boldsymbol{p}_h) .
\end{align*}
By \eqref{norm_sigma}, the projected gradient term controls the broken \(H^1\) seminorm of the virtual tensor, so no additional stabilization is needed for this variable. The term \((\operatorname{div}\boldsymbol r_h,\operatorname{div}\boldsymbol s)\) supplies the divergence control required by \eqref{eq_rt_equivalence} and introduces no consistency error, since the exact multiplier satisfies \(\boldsymbol r=0\) by Theorem~\ref{thm:gStokes_continue_wellposed}.

Let $Q_h^{\mathrm{nt}}:L^2(\Omega;\mathbb T)\to\Vhnt$ denote the $L^2$-orthogonal projection.

\begin{lemma}
Let $k\ge0$. For $\left( \boldsymbol{\tau}, \boldsymbol{s}  \right) \in \VhvemS   \times \VhRT/\mathbb{R}^3 $ satisfying $ b_h \left( \boldsymbol{\tau}, \boldsymbol{s}; \boldsymbol{q} \right) = 0 $ for all $\boldsymbol{q} \in \Vhnt $, we have 
\begin{equation}\label{coer}
 \| \gradvem \bs \tau \|^2 + \interleave \boldsymbol{s} \interleave_{1, h}^2  \lesssim a_h\left(\boldsymbol{\tau}, \boldsymbol{s}; \boldsymbol{\tau}, \boldsymbol{s} \right).
\end{equation} 
\end{lemma}
\begin{proof}
By the assumption \(b_h(\boldsymbol\tau,\boldsymbol s;\boldsymbol q)=0\) for all \(\boldsymbol q\in\Vhnt\), we have
\[
\|\dev\grad_w\bs s\|=\|\Qnt\curlvem\bs\tau\|
\lesssim\|\gradvem\bs\tau\|.
\]
Combining this bound with \eqref{eq_rt_equivalence} proves \eqref{coer}.
\end{proof}

Let $\Ivem:H_0^1(\Omega;\bbS)\to\VhvemS$, $\Int: H^1(\Omega; \mathbb{T})\cap H_0(\div,\Omega;\mathbb{T}) \to \Vhnt $ and $I_h^{\mathrm{RT}}:H^1(\Omega;\bbR^3)\to
\VhRT$ be the canonical interpolation operators defined by the DoFs in \eqref{dof_vem}, \eqref{dof_nt} and \eqref{dof_hdiv}, respectively.
For $0\leq m\leq k+1$, let $\bs\tau\in H^{m+1}(\Omega;\mathbb S)\cap H_0^1(\Omega;\mathbb S)$. For $1\leq m\leq k+1$, let $\bs s\in H^m(\Omega;\mathbb R^3)$ and $\bs p\in H^m(\Omega;\mathbb T)\cap H_0(\div,\Omega;\mathbb T)$. Standard approximation and scaling estimates for the nonconforming virtual element \cite{ChenHuang2020,HuangTang2025}, tangential-normal element \cite{ChenHuangZhao2026}, and Raviart--Thomas element \cite{RaviartThomas1977,Nedelec1980,ChenHuang2022b} interpolants give
\begin{align}
\notag
\|\grad\bs\tau-\gradvem(\Ivem\bs\tau)\|
\lesssim h^m|\bs\tau|_{m+1}, \\
  \label{interpolation_error_p}
\|\bs p-\Int\bs p\|+h|\bs p-\Int\bs p|_{1,h}
\lesssim h^m|\bs p|_m, \\
\label{interpolation_error_s}
\|\bs s-I_h^{\mathrm{RT}}\bs s\|+h\interl{\bs s-I_h^{\mathrm{RT}}\bs s}
\lesssim h^m|\bs s|_m.
\end{align}
The local VE interpolant $I_T^{\mathrm{VE}}$ also satisfies \cite[Section~3.3]{HuangTang2025}, for $0\le s\le k+1$, 
\begin{equation}\label{eq:VE_local_interpolation}
\|\bs\tau-I_T^{\mathrm{VE}}\bs\tau\|_T
+h_T|\bs\tau-I_T^{\mathrm{VE}}\bs\tau|_{1,T}
\lesssim h_T^{s+1}|\bs\tau|_{s+1,T}, \quad\forall\,\bs\tau\in H^{s+1}(T;\mathbb S).
\end{equation}

\begin{lemma}\label{lem:commuting_properties}
For every integer $k\ge0$, $\bs\tau\in H_0^1(\Omega;\bbS)$, and $\bs s\in H^1(\Omega;\bbR^3)$, the following identities hold:
\begin{subequations}
\begin{align}
\gradvem(\Ivem \bs \tau) & =  Q_{k,h}(\grad \bs \tau), \label{com_grad} \\
\curlvem(\Ivem \bs \tau) &= Q_{k,h} (\curl \bs \tau),  \label{com_curl} \\
\dg_w(I_h^{\mathrm{RT}} \bs s) &= \Qnt (\dg \bs s). \label{com_weak_grad}
\end{align}
\end{subequations}
\end{lemma}
\begin{proof}
Identities \eqref{com_grad}--\eqref{com_curl} follow from the commuting property of the nonconforming VE interpolant \cite[Section~3.3]{HuangTang2025}, applied componentwise. For \eqref{com_weak_grad}, the RT moments \eqref{dof_hdiv} and \eqref{eq_weak_devgrad} give, for every $\bs q_h\in\Vhnt$,
\[
\begin{aligned}
(\dev\grad_w(I_h^{\mathrm{RT}}\bs s),\bs q_h)
&=-(\div_h\bs q_h,\bs s)
+(\jump{\bs n_F^{\intercal}\bs q_h\bs n_F},\bs s\cdot\bs n_F)_{\mathcal F_h}\\
&=(\dev\grad\bs s,\bs q_h).
\end{aligned}
\]
The last equality uses elementwise integration by parts, tangential-normal continuity and its homogeneous boundary condition, and the tracelessness of $\bs q_h$. This proves \eqref{com_weak_grad}.
\end{proof}

\begin{lemma}
For every fixed integer $k\ge0$, the following \textup{inf-sup} condition holds:
\begin{equation}\label{discrete_inf-sup}
\|\bs q_h\|
\lesssim
\sup_{\bs\tau_h\in\VhvemS,\ \bs s_h\in \VhRT/\mathbb{R}^3 }
\frac{ b_h(\bs\tau_h,\bs s_h;\bs q_h)}{\| \gradvem \bs \tau_h \| + \interl{\bs s_h} },
\qquad   \bs q_h\in\Vhnt.
\end{equation}
By the norm equivalence \eqref{eq_rt_equivalence}, $\interl{\bs s_h}$ can be replaced by $\|\bs s_h\|_{R,h}$.
\end{lemma}

\begin{proof}
The Helmholtz decomposition for traceless tensors \cite[Lemma~2]{GuHuang2026} gives 
\begin{equation}\label{eq:qhhelmholtz_decomposition}
\bs q_h = \curl \bs \tau + \dev\grad\bs s
\end{equation}
with $\bs \tau \in H_0^1(\Omega;\bbS)$ and $\bs s \in H^1(\Omega;\bbR^3)/\text{RT}$ satisfying
\begin{equation*}
\| \bs \tau \|_1 + \| \bs s \|_1 \lesssim \| \bs q_h \|.
\end{equation*}
Choose
$\bs\tau_h=\Ivem\bs\tau$ and $\bs s_h=I_h^{\mathrm{RT}}\bs s-\frac1{|\Omega|}\int_\Omega I_h^{\mathrm{RT}}\bs s\,dx$, then $\bs s_h\in\VhRT/\mathbb R^3$.  Applying $Q_h^{\mathrm{nt}}$ to the decomposition \eqref{eq:qhhelmholtz_decomposition} and using \eqref{com_curl}--\eqref{com_weak_grad} gives
\[
\bs q_h=Q_h^{\mathrm{nt}}\curlvem\bs\tau_h+\dev\grad_w\bs s_h.
\]
Consequently, $b_h(\bs\tau_h,\bs s_h;\bs q_h)=\|\bs q_h\|^2$.
The commuting identities, the $L^2$ stability of the projections, and \eqref{interpolation_error_s} also give
\[
\|\gradvem\bs\tau_h\|+\interl{\bs s_h}
\lesssim \|\bs\tau\|_1+\|\bs s\|_1
\lesssim \|\bs q_h\|.
\]
\end{proof}


To obtain the optimal \(L^2\) convergence rate for \(k=1\), we replace the virtual element test function in the load term by a piecewise-linear reconstruction. For a unified treatment of all \(k\ge0\), we define the reconstructed test function elementwise for \(\boldsymbol\tau\in H^1(\mathcal T_h;\mathbb S)\) by

\begin{equation}\label{eq:load_reconstruction}
\widetilde{\bs\tau}|_T=
\begin{cases}
I_T^{\mathrm{VE}}\bs\tau,&k=0,\\
Q_{0,T}\bs\tau+(Q_{0,T}\nabla\bs\tau)(\bs x-\bs x_T),&k=1,\\
Q_{k-1,T}\bs\tau,&k\ge2.
\end{cases}
\end{equation}
For $k=0$, $I_T^{\mathrm{VE}}$ is the CR interpolant, so $\widetilde{\bs\tau}_h=\bs\tau_h$ for discrete test tensors. For \(k=1\), the reconstruction is the elliptic projection onto \(\mathbb P_1(T;\mathbb S)\) that preserves the element mean, and it is computable from the DoFs~\eqref{dof_vem}.

The moment preservation property, the idempotence of the interpolation operator, linear polynomial reproduction, and \eqref{eq:VE_local_interpolation} give, for all $k\ge0$,
\begin{subequations}\label{eq:load_properties}
\begin{equation}\label{eq:load_moments}
Q_{k-1,T}\widetilde{\bs\tau}=Q_{k-1,T}\bs\tau,
\qquad
\widetilde{\Ivem\bs\tau}=\widetilde{\bs\tau},
\end{equation}
where the second identity holds for \(\boldsymbol\tau\in H_0^1(\Omega;\mathbb S)\), the domain of the global interpolation operator \(I_h^{\mathrm{VE}}\), and
\begin{equation}\label{eq:load_approximation}
\|\bs\tau-\widetilde{\bs\tau}\|_T\lesssim h_T^2|\bs\tau|_{2,T},
\qquad\bs\tau\in H^2(T;\mathbb S).
\end{equation}
\end{subequations}
The definition of the reconstruction, the \(L^2\) stability of the projections, and \eqref{eq:VE_local_interpolation} with \(s=0\) yield
\begin{equation}\label{eq:load_stability}
\|\widetilde{\bs\tau}\|\lesssim\|\bs\tau\|+h|\bs\tau|_{1,h},
\qquad\bs\tau\in H^1(\mathcal T_h;\mathbb S).
\end{equation}
We use $(\bs g_h,\widetilde{\bs\tau})$ as the discrete load.
The Cauchy--Schwarz inequality, $\|\curlvem\bs\tau\|\lesssim\|\gradvem\bs\tau\|$, and the norm equivalence \eqref{eq_rt_equivalence} give the discrete continuity estimates
\begin{subequations}\label{eq:discrete_continuity}
\begin{align}
|a_h(\bs\sigma,\bs r;\bs\tau,\bs s)|
&\lesssim (\|\gradvem\bs\sigma\|+\interl{\bs r})
(\|\gradvem\bs\tau\|+\interl{\bs s}),\\
|b_h(\bs\tau,\bs s;\bs p)|
&\lesssim (\|\gradvem\bs\tau\|+\interl{\bs s})\|\bs p\|,
\end{align}
\end{subequations}
for $(\bs\sigma,\bs r),(\bs\tau,\bs s)\in\VhvemS\times(\VhRT/\mathbb R^3)$ and $\bs p\in\Vhnt$.
\begin{theorem}\label{thm:gStokes_discrete_wellposed}
For every integer $k\ge0$, the distributional finite element method \eqref{discrete_generalized_stokes} is well-posed.
\end{theorem}
\begin{proof}
The load stability \eqref{eq:load_stability}, the broken Poincar\'e inequality, and \eqref{norm_sigma} give
\[
|(\bs g_h,\widetilde{\bs\tau})|\lesssim\|\bs g_h\|\|\gradvem\bs\tau\|,
\qquad k\ge0.
\]
The discrete continuity estimates \eqref{eq:discrete_continuity}, coercivity on the kernel of bilinear form $b_h$ \eqref{coer}, and inf-sup condition \eqref{discrete_inf-sup} then imply well-posedness by the Babuška--Brezzi theory \cite{BoffiBrezziFortin2013}.
\end{proof}

\subsection{Optimal energy and $L^2$ error estimates}
\begin{lemma}\label{lem_err1}
Let $k\ge0$ and $(\bs\sigma,\bs p,0)$ solve \eqref{continue_generalized_stokes}, with $\bs\sigma\in H^{k+2}(\Omega;\mathbb S)$ and $\bs p\in H^{k+1}(\Omega;\mathbb T)$. For $(\bs\tau,\bs s)\in\VhvemS\times(\VhRT/\mathbb R^3)$, define
\[
\mathcal E_h(\bs\tau,\bs s)
:=a_h(\Ivem\bs\sigma,0;\bs\tau,\bs s)
+b_h(\bs\tau,\bs s;\Int\bs p)-(\bs g,\widetilde{\bs\tau}).
\]
Then
\begin{equation}\label{eq_err_estimate1}
|\mathcal E_h(\bs\tau,\bs s)|
\lesssim h^{k+1}(|\bs\sigma|_{k+2}+|\bs p|_{k+1})
\bigl(\|\gradvem\bs\tau\|+\interl{\bs s}\bigr).
\end{equation}
\end{lemma}
\begin{proof}
By Theorem~\ref{thm:gStokes_continue_wellposed}, $\bs p \in H_0(\div,\Omega;\mathbb T)\cap\ker(\div)$. Thus definition~\eqref{eq_weak_devgrad} of $\dev\grad_w$ and elementwise integration by parts yield
\[
(\Int\bs p,\dev\grad_w\bs s)
=(\Int\bs p-\bs p,\dev\grad_h\bs s)
-(\Pi_F(\Int\bs p-\bs p)\bs n_F,\jump{\Pi_F\bs s})_{\innerF}.
\]
The trace inequality and interpolation estimate \eqref{interpolation_error_p} give
\begin{equation}\label{eq_err5}
|(\Int\bs p,\dev\grad_w\bs s)|
\lesssim h^{k+1}|\bs p|_{k+1}\interl{\bs s}.
\end{equation}

For $k=0$, $\Ivem$ is the CR interpolant, $\widetilde{\bs\tau}=\bs\tau$, and $\gradvem\bs\tau=\grad_h\bs\tau$. The commuting identity \eqref{com_grad} and elementwise integration by parts give
\[
\begin{aligned}
\mathcal E_h(\bs\tau,\bs s)
={}&(\grad\bs\sigma\,\bs n+\bs p\times\bs n,\bs\tau)_{\partial\mathcal T_h}
+(\Int\bs p-\bs p,\curl_h\bs\tau)\\
&+(\Int\bs p,\dev\grad_w\bs s).
\end{aligned}
\]
Applying the continuity of face means and the trace inequality to the first term, the interpolation estimate for \(I_h^{\mathrm{nt}}\boldsymbol p\) to the second term, and \eqref{eq_err5} to the third term yields
\[
|\mathcal E_h(\bs\tau,\bs s)|
\lesssim h(|\bs\sigma|_2+|\bs p|_1)
\bigl(|\bs\tau|_{1,h}+\interl{\bs s}\bigr),
\]
which is \eqref{eq_err_estimate1} for $k=0$.

For $k\ge1$, elementwise integration by parts, the commuting identity~\eqref{com_grad}, and the cell-moment preservation of $\Int$ yield
\begin{equation}\label{eq_err2}
\begin{aligned}
(\grad\bs\sigma,\grad_h\bs\tau_T)
&=(\gradvem\Ivem\bs\sigma,\gradvem\bs\tau)
+(Q_{k,h}\grad\bs\sigma\,\bs n,\bs\tau_T-\bs\tau_F)_{\partial\mathcal T_h},\\
(\bs p,\curl_h\bs\tau_T)
&=(\Int\bs p,\curlvem\bs\tau)
+((\Int\bs p)\times\bs n,\bs\tau_T-\bs\tau_F)_{\partial\mathcal T_h}.
\end{aligned}
\end{equation}
For simplicity, define
\[
R_h(\bs\tau):=
\bigl((Q_{k,h}-I)\grad\bs\sigma\,\bs n
+(\Int\bs p-\bs p)\times\bs n,\bs\tau_T-\bs\tau_F\bigr)_{\partial\mathcal T_h}.
\]
Using $\bs g=-\Delta\bs\sigma+\sym\curl\bs p$, \eqref{eq_err2}, and the cancellation of the exact-flux face terms, we obtain
\begin{equation}\label{eq_err4}
(\bs g,\bs\tau_T)
=(\gradvem\Ivem\bs\sigma,\gradvem\bs\tau)
+(\Int\bs p,\curlvem\bs\tau)+R_h(\bs\tau).
\end{equation}
Consequently,
\[
\mathcal E_h(\bs\tau,\bs s)
=-R_h(\bs\tau)+(\Int\bs p,\dev\grad_w\bs s)
+(\bs g,\bs\tau_T-\widetilde{\bs\tau}).
\]
The interpolation estimates and \eqref{eq:vem_trace_estimate} imply
\begin{equation}\label{eq_err6}
|R_h(\bs\tau)|\lesssim
h^{k+1}(|\bs\sigma|_{k+2}+|\bs p|_{k+1})\|\gradvem\bs\tau\|.
\end{equation}
For $k\ge2$, the load correction term $ (\bs g,\bs\tau_T-\widetilde{\bs\tau}) $ vanishes. For $k=1$, \eqref{eq:load_reconstruction}--\eqref{eq:load_moments} and the elementwise Poincar\'e inequality give
\[
\begin{aligned}
|(\bs g,\bs\tau_T-\widetilde{\bs\tau})|
&=|(\bs g-Q_{0,h}\bs g,\bs\tau_T-\widetilde{\bs\tau})|\\
&\lesssim h^2|\bs g|_1\|\gradvem\bs\tau\|
\lesssim h^2(|\bs\sigma|_3+|\bs p|_2)\|\gradvem\bs\tau\|.
\end{aligned}
\]
Combining these bounds proves \eqref{eq_err_estimate1}.
\end{proof}
\begin{theorem}\label{thm_H1bound}
Let $k\ge0$, and let $(\bs\sigma,\bs p,0)$ and $(\bs\sigma_h,\bs p_h,\bs r_h)$ solve \eqref{continue_generalized_stokes} and \eqref{discrete_generalized_stokes}, respectively. Assume $\bs\sigma\in H^{k+2}(\Omega;\mathbb S)$ and $\bs p\in H^{k+1}(\Omega;\mathbb T)$. Then
\begin{equation}\label{err_stokes_H1}
\begin{aligned} 
& \| Q_{k,h}( \grad\bs \sigma) - \gradvem \boldsymbol{\sigma}_h \| + \| \boldsymbol{p} - \boldsymbol{p}_h \| + \interleave\boldsymbol{r}_h\interleave_{1, h} \\
&\qquad\qquad\qquad\qquad\qquad\qquad\lesssim h^{k+1}(| \boldsymbol{\sigma} |_{k+2} +| \boldsymbol{p} |_{k+1}) + \| \bs g - \bs g_h \|.
\end{aligned}
\end{equation}
\end{theorem}
\begin{proof}
Set $\bs e_\sigma:=\bs\sigma_h-\Ivem\bs\sigma$ and $\bs e_p:=\bs p_h-\Int\bs p$. Since $\curl\bs\sigma=0$, the commuting identity \eqref{com_curl} gives $\curlvem\Ivem\bs\sigma=0$. Subtracting the interpolated equations from \eqref{discrete_generalized_stokes} yields, for all $(\bs\tau,\bs s)\in\VhvemS\times(\VhRT/\mathbb R^3)$ and $\bs q\in\Vhnt$,
\begin{equation}\label{eq:error_equations}
\begin{aligned}
a_h(\bs e_\sigma,\bs r_h;\bs\tau,\bs s)+b_h(\bs\tau,\bs s;\bs e_p)
&=(\bs g_h-\bs g,\widetilde{\bs\tau})-\mathcal E_h(\bs\tau,\bs s),\\
b_h(\bs e_\sigma,\bs r_h;\bs q)&=0.
\end{aligned}
\end{equation}
In particular, the second equation gives the constraint used below:
\begin{equation}\label{eq_err0}
\Qnt\curlvem\bs e_\sigma+\dev\grad_w\bs r_h=0.
\end{equation}
The coercivity estimate~\eqref{coer}, the norm equivalence~\eqref{eq_rt_equivalence}, and Lemma~\ref{lem_err1} give
\[
\begin{aligned}
&\|\operatorname{grad}_h^\Pi\boldsymbol e_\sigma\|^2+\|\boldsymbol r_h\|_{R,h}^2
\lesssim
a_h(\boldsymbol e_\sigma,\boldsymbol r_h;\boldsymbol e_\sigma,\boldsymbol r_h) =
(\boldsymbol g_h-\boldsymbol g,\widetilde{\boldsymbol e_\sigma})
-\mathcal E_h(\boldsymbol e_\sigma,\boldsymbol r_h)\\
&\qquad\qquad\lesssim
\bigl[h^{k+1}(|\boldsymbol\sigma|_{k+2}+|\boldsymbol p|_{k+1})
+\|\boldsymbol g-\boldsymbol g_h\|\bigr]
\bigl(\|\operatorname{grad}_h^\Pi\boldsymbol e_\sigma\|
+\|\boldsymbol r_h\|_{R,h}\bigr).
\end{aligned}
\]
In addition, the inf-sup condition~\eqref{discrete_inf-sup}, the error equations~\eqref{eq:error_equations}, and the discrete continuity~\eqref{eq:discrete_continuity} yield
\[
\begin{aligned}
\|\boldsymbol e_p\|
&\lesssim
h^{k+1}(|\boldsymbol\sigma|_{k+2}+|\boldsymbol p|_{k+1})
+\|\boldsymbol g-\boldsymbol g_h\|
+\|\operatorname{grad}_h^\Pi\boldsymbol e_\sigma\|
+\|\boldsymbol r_h\|_{R,h}\\
&\lesssim
h^{k+1}(|\boldsymbol\sigma|_{k+2}+|\boldsymbol p|_{k+1})
+\|\boldsymbol g-\boldsymbol g_h\|.
\end{aligned}
\]
Combining these estimates, we have 
\begin{equation}\label{eq_err_H1}
\begin{aligned}
\|\gradvem\bs e_\sigma\|+\|\bs r_h\|_{R,h}+\|\bs e_p\|
\lesssim h^{k+1}(|\bs\sigma|_{k+2}+|\bs p|_{k+1})+\|\bs g-\bs g_h\|.
\end{aligned}
\end{equation}
Finally, the interpolation estimate~\eqref{interpolation_error_p} for \(\boldsymbol p\), the commuting identity~\eqref{com_grad}, and the Raviart–Thomas norm equivalence~\eqref{eq_rt_equivalence} yield \eqref{err_stokes_H1}.
\end{proof}

In particular, combining \eqref{err_stokes_H1} with the projection estimate of $Q_{k,h}$ gives
\begin{equation}\label{eq:stokes_full_gradient}
\|\grad\bs\sigma-\gradvem\bs\sigma_h\|+\|\bs p-\bs p_h\|+\interl{\bs r_h}
\lesssim h^{k+1}(|\bs\sigma|_{k+2}+|\bs p|_{k+1})+\|\bs g-\bs g_h\|.
\end{equation}

We then give the optimal $L^2$ estimate of $\bs \sigma-\bs \sigma_h$ through a duality argument. Let $\bs e:=\Ivem\bs\sigma-\bs\sigma_h\in\VhvemS$. For $k\ge1$, define
\[
\bs e_T:=Q_{k-1,h}\bs e=\bs\sigma_T-\bs\sigma_{h,T},\qquad
\bs e_F:=Q_{k,\mcal F_h}\bs e=\bs\sigma_F-\bs\sigma_{h,F}.
\]
Let $\hat u\in H_0^3(\Omega)$ solve
\begin{equation}\label{dual_problem}
-\Delta^3\hat u=\div\div\bs e_T\quad\text{in the distributional sense}.
\end{equation}
We assume the dual regularity estimates
\begin{subequations}\label{dual_regularity}
\begin{equation}\label{dual_regu}
\|\hat u\|_4\lesssim\|\bs e_T\|,
\end{equation}
\begin{equation}\label{traceless_regu}
\|\bs q\|_1\lesssim\|\div\bs q\|+\|\sym\curl\bs q\|,\quad
\bs q\in H_0(\div,\Omega;\mathbb T)\cap H(\sym\curl,\Omega;\mathbb T).
\end{equation}
\end{subequations}
Here \(
H(\sym\curl,\Omega;\mathbb T)
:=\{\bs q\in L^2(\Omega;\mathbb T):\sym\curl\bs q\in L^2(\Omega;\mathbb S)\}.
\)
We refer to \cite{MazyaRossmann2010,KozlovMazyaRossmann2001} for regularity results for polyharmonic equations; related estimates for vector fields on convex domains can be found in \cite[Section~3.5]{GiraultRaviart1986} and \cite[Corollary~5.2]{Mitrea2001}. 
In particular, we do not claim \eqref{dual_regularity} on arbitrary polyhedral domains.

Applying Theorem~\ref{thm:gStokes_continue_wellposed} with the load $\bs e_T$, and using the exactness of \eqref{hessiancomplex} and the regularity assumptions \eqref{dual_regularity}, yields
\(
\hat{\bs\sigma}=\nabla^2\hat u\in H_0^1(\Omega;\mathbb S)\cap H^2(\Omega;\mathbb S)\) and
\(\hat{\bs p}\in H_0(\div,\Omega;\mathbb T)\cap H^1(\Omega;\mathbb T)\)
such that
\begin{equation}\label{dual_replace}
-\Delta\hat{\bs\sigma}+\sym\curl\hat{\bs p}=\bs e_T,\qquad
\curl\hat{\bs\sigma}=0,\qquad\div\hat{\bs p}=0,
\end{equation}
and
\begin{equation}\label{dual_estimate}
\|\hat{\bs\sigma}\|_2+\|\hat{\bs p}\|_1\lesssim\|\bs e_T\|.
\end{equation}

\begin{theorem}\label{thm:stokes_L2_error}
Let $k\ge1$. Under the assumptions of Theorem~\ref{thm_H1bound} and the dual regularity assumptions \eqref{dual_regularity}, we have 
\begin{equation}\label{err_stokes_L2}
\begin{aligned}
\|\bs\sigma-\bs\sigma_h\|
\lesssim h^{k+2}(|\bs\sigma|_{k+2}+|\bs p|_{k+1}) +h\|\bs g-\bs g_h\|+\|\bs g-\bs g_h\|_{-1}.
\end{aligned}
\end{equation}
\end{theorem}
\begin{proof}
We separately estimate the two terms in
\[
\|\bs e_T\|^2=(\bs e_T,-\Delta\hat{\bs\sigma})
+(\bs e_T,\sym\curl\hat{\bs p}).
\]

For the second term, elementwise integration by parts and the second identity in \eqref{eq_err2}, applied with $(\bs\tau,\bs p)=(\bs e,\hat{\bs p})$, give
\[
(\bs e_T,\sym\curl\hat{\bs p})
=(\curlvem\bs e,\Int\hat{\bs p})+\big((\Int\hat{\bs p}-\hat{\bs p})\times\bs n,\bs e_T-\bs e_F\big)_{\partial\mathcal T_h}.
\]
Moreover, \eqref{eq_err0} and $\bs e=-\bs e_\sigma$ give
\[
(\curlvem\bs e,\Int\hat{\bs p})= (\Qnt\curlvem\bs e,\Int\hat{\bs p})= (\dg_w\bs r_h,\Int\hat{\bs p}).
\]
Since \(\hat{\boldsymbol p}\in H_0(\operatorname{div},\Omega;\mathbb T)\cap\ker(\operatorname{div})\), the argument used to prove \eqref{eq_err5}, together with the interpolation estimate~\eqref{interpolation_error_p} for \(m=1\), yields
\[
\bigl|(\curlvem\bs e,\Int\hat{\bs p})\bigr|= \bigl|(\dg_w\bs r_h,\Int\hat{\bs p})\bigr|
\lesssim h\|\hat{\bs p}\|_1\|\bs r_h\|_{R,h}.
\]
On the other hand,  the interpolation estimate~\eqref{interpolation_error_p} with \(m=1\), and the virtual element trace estimate~\eqref{eq:vem_trace_estimate} yield
\[
\left|
\bigl((I_h^{\mathrm{nt}}\hat{\boldsymbol p}-\hat{\boldsymbol p})\times\boldsymbol n,
\boldsymbol e_T-\boldsymbol e_F\bigr)_{\partial\mathcal T_h}
\right|
\lesssim h\|\hat{\boldsymbol p}\|_1
\|\operatorname{grad}_h^\Pi\boldsymbol e\|.
\]

Consequently, these two bounds, together with the energy error estimate~\eqref{eq_err_H1} and the dual regularity estimate~\eqref{dual_estimate}, imply
\begin{equation}\label{pre:dualestimate1}
|(\boldsymbol e_T,\operatorname{sym}\operatorname{curl}\hat{\boldsymbol p})|
\lesssim
\|\boldsymbol e_T\|
\bigl[
h^{k+2}(|\boldsymbol\sigma|_{k+2}+|\boldsymbol p|_{k+1})
+h\|\boldsymbol g-\boldsymbol g_h\|
\bigr].
\end{equation}

For the first term, applying the first identity in \eqref{eq_err2} with \(\boldsymbol\sigma=\hat{\boldsymbol\sigma}\) and \(\boldsymbol\tau=\boldsymbol e\) gives
\begin{equation}
\begin{aligned}
(\bs e_T,-\Delta\hat{\bs\sigma})
= ( \gradvem \bs e, \gradvem \Ivem \hat{\bs\sigma})
+( (Q_{k,h} - I ) \nabla \hat{ \bs\sigma} \, \bs n,
\bs e_T-\bs e_F)_{\partial\mathcal T_h}. \nonumber
\end{aligned}
\end{equation}
Denote $A_h:=(\operatorname{grad}_h^\Pi\boldsymbol e,
\operatorname{grad}_h^\Pi I_h^{\mathrm{VE}}\hat{\boldsymbol\sigma}).$ 
For any $k\geq0$,
since \(\operatorname{curl}\hat{\boldsymbol\sigma}=0\), we test the discrete equation~\eqref{discrete_generalized_stokes1} with \((I_h^{\mathrm{VE}}\hat{\boldsymbol\sigma},0)\). The commuting identities~\eqref{com_grad}-\eqref{com_curl} and the reconstruction identity~\eqref{eq:load_moments}, together with \(\boldsymbol e=I_h^{\mathrm{VE}}\boldsymbol\sigma-\boldsymbol\sigma_h\), give
\begin{equation}\label{eq:dual_testing_identity}
\begin{aligned}
A_h &=(Q_{k,h}\nabla\boldsymbol\sigma,Q_{k,h}\nabla\hat{\boldsymbol\sigma}) 
-(\boldsymbol g_h,\widetilde{\hat{\boldsymbol\sigma}})\\
&= -((I-Q_{k,h})\nabla\boldsymbol\sigma,
(I-Q_{k,h})\nabla\hat{\boldsymbol\sigma})
+(\boldsymbol g,\hat{\boldsymbol\sigma}-\widetilde{\hat{\boldsymbol\sigma}})
+(\boldsymbol g-\boldsymbol g_h,\widetilde{\hat{\boldsymbol\sigma}}).
\end{aligned}
\end{equation}
The second equality follows from \eqref{continue_generalized_stokes1} and the $L^2$-orthogonality of $Q_{k,h}$.
The Cauchy–Schwarz inequality and the polynomial projection approximation estimates yield
\[
\begin{aligned}
\big|((I-Q_{k,h})\nabla\boldsymbol\sigma,
(I-Q_{k,h})\nabla\hat{\boldsymbol\sigma})\big|
\lesssim h^{k+2}|\boldsymbol\sigma|_{k+2}|\hat{\boldsymbol\sigma}|_2.
\end{aligned}
\]
The moment-preserving property~\eqref{eq:load_moments}, the polynomial projection approximation estimate, and the reconstruction approximation estimate~\eqref{eq:load_approximation} give
\begin{equation}\label{eq_L2_tmp}
\begin{aligned}
\big|(\bs g,\hat{\bs\sigma}-\widetilde{\hat{\bs\sigma}})\big|
&=\big|(\bs g-Q_{k-1,h}\bs g,\hat{\bs\sigma}-\widetilde{\hat{\bs\sigma}})\big| \lesssim h^{k+2}(|\bs\sigma|_{k+2}+|\bs p|_{k+1})|\hat{\bs\sigma}|_2, \\
|(\boldsymbol g-\boldsymbol g_h,\widetilde{\hat{\boldsymbol\sigma}})|
&=\big|(\boldsymbol g-\boldsymbol g_h,\widetilde{\hat{\boldsymbol\sigma}}-\hat{\boldsymbol\sigma})
+(\boldsymbol g-\boldsymbol g_h,\hat{\boldsymbol\sigma})\big|\\
&\lesssim h^2\|\boldsymbol g-\boldsymbol g_h\|\,|\hat{\boldsymbol\sigma}|_2
+\|\boldsymbol g-\boldsymbol g_h\|_{-1}\|\hat{\boldsymbol\sigma}\|_1.
\end{aligned}
\end{equation}
For the boundary term, the virtual element trace estimate~\eqref{eq:vem_trace_estimate} and the energy error estimate~\eqref{err_stokes_H1} yield
\[
\begin{aligned}
\big|
\bigl((Q_{k,h}-I)\nabla\hat{\boldsymbol\sigma}\,\boldsymbol n,
\boldsymbol e_T-\boldsymbol e_F\bigr)_{\partial\mathcal T_h}
\big|
&\lesssim
h|\hat{\boldsymbol\sigma}|_2\|\operatorname{grad}_h^\Pi\boldsymbol e\|\\
&\lesssim
|\hat{\boldsymbol\sigma}|_2
\bigl[h^{k+2}(|\boldsymbol\sigma|_{k+2}+|\boldsymbol p|_{k+1})
+h\|\boldsymbol g-\boldsymbol g_h\|\bigr].
\end{aligned}
\]
Combining these bounds with the preceding estimates for \(A_h\), the dual regularity estimate~\eqref{dual_estimate}, and \(h\lesssim1\), we obtain
\begin{equation}\label{pre:dualestimate2}
|(\boldsymbol e_T,-\Delta\hat{\boldsymbol\sigma})|
\lesssim
\|\boldsymbol e_T\|
\bigl[
h^{k+2}(|\boldsymbol\sigma|_{k+2}+|\boldsymbol p|_{k+1})
+h\|\boldsymbol g-\boldsymbol g_h\|
+\|\boldsymbol g-\boldsymbol g_h\|_{-1}
\bigr].
\end{equation}
Combining the estimates~\eqref{pre:dualestimate1} and \eqref{pre:dualestimate2}, we have
\begin{equation}
\|\boldsymbol e_T\|
\lesssim
h^{k+2}(|\boldsymbol\sigma|_{k+2}+|\boldsymbol p|_{k+1})
+h\|\boldsymbol g-\boldsymbol g_h\|
+\|\boldsymbol g-\boldsymbol g_h\|_{-1}. \nonumber
\end{equation}
Since \(k\ge1\), the projection \(Q_{k-1,T}\) preserves the element mean. The local Poincaré inequality and the projected gradient norm equivalence~\eqref{norm_sigma} therefore give
\[
\|\boldsymbol e\|
\lesssim\|Q_{k-1,h}\boldsymbol e\|+h|\boldsymbol e|_{1,h}
\lesssim\|\boldsymbol e_T\|+h\|\operatorname{grad}_h^\Pi\boldsymbol e\|.
\]
Finally, using \(\boldsymbol e=I_h^{\mathrm{VE}}\boldsymbol\sigma-\boldsymbol\sigma_h\), the interpolation estimate~\eqref{eq:VE_local_interpolation}, and the energy error estimate~\eqref{eq_err_H1}, we obtain
\[
\begin{aligned}
	\|\boldsymbol\sigma-\boldsymbol\sigma_h\|
	&\le\|\boldsymbol\sigma-I_h^{\mathrm{VE}}\boldsymbol\sigma\|+\|\boldsymbol e\|\\
	&\lesssim
	h^{k+2}(|\boldsymbol\sigma|_{k+2}+|\boldsymbol p|_{k+1})
	+h\|\boldsymbol g-\boldsymbol g_h\|
	+\|\boldsymbol g-\boldsymbol g_h\|_{-1},
\end{aligned}
\]
which proves \eqref{err_stokes_L2}.
\end{proof}

From Theorem~\ref{thm:stokes_L2_error}, the computable element moments inherit the same error bound by $L^2$-projection stability:
\begin{equation}\label{compulable_L2_bound}
\|\bs\sigma_T-\bs\sigma_{h,T}\|
=\|Q_{k-1,h}(\bs\sigma-\bs\sigma_h)\|
\le\|\bs\sigma-\bs\sigma_h\|,\qquad k\ge1.
\end{equation}

\begin{remark}
If the discrete load is evaluated exactly, i.e., $(\bs g_h,\widetilde{\bs\tau})=(\bs g,\widetilde{\bs\tau})$ for every discrete test tensor $\bs\tau$, the discrete solution coincides with that obtained by taking $\bs g_h=\bs g$, hence Theorem~\ref{thm:stokes_L2_error} gives the optimal $O(h^{k+2})$ error for the full virtual tensor under its regularity assumptions. For $k=1$, replacing the reconstructed load $(\bs g,\widetilde{\bs\tau})$ by $(\bs g,\bs\tau_T)$ gives only an $O(h^2)$ bound in the preceding analysis.
\end{remark}

We also have the following $ L^2 $ estimate for $ k = 0$.
\begin{proposition}\label{prop:k0_errors}
Under the assumptions of Theorem~\ref{thm_H1bound} with $k=0$, assume that there exist $\hat{\bs\sigma}\in H_0^1(\Omega;\mathbb S)\cap H^2(\Omega;\mathbb S)$ and $\hat{\bs p}\in H_0(\div,\Omega;\mathbb T)\cap H^1(\Omega;\mathbb T)$ satisfying
\begin{equation}\label{eq:k0_dual}
\begin{gathered}
-\Delta\hat{\bs\sigma}+\sym\curl\hat{\bs p}=\bs\sigma-\bs\sigma_h,
\qquad \curl\hat{\bs\sigma}=0,\qquad \div\hat{\bs p}=0,\\
\|\hat{\bs\sigma}\|_2+\|\hat{\bs p}\|_1\lesssim\|\bs\sigma-\bs\sigma_h\|,
\end{gathered}
\end{equation}
then the $L^2$ error satisfies
\begin{equation}\label{eq:k0_L2}
\|\bs\sigma-\bs\sigma_h\|
\lesssim h^2(|\bs\sigma|_2+|\bs p|_1)
+h\|\bs g-\bs g_h\|+\|\bs g-\bs g_h\|_{-1}.
\end{equation}
\end{proposition}
\begin{proof}
Testing the dual equation~\eqref{eq:k0_dual} with \(\boldsymbol\sigma-\boldsymbol\sigma_h\) elementwise gives
\begin{equation}\label{eq:prop_tmp}
\begin{aligned}
\|\boldsymbol\sigma-\boldsymbol\sigma_h\|^2
={}&(\nabla_h(\boldsymbol\sigma-\boldsymbol\sigma_h),
\nabla_h I_h^{\mathrm{VE}}\hat{\boldsymbol\sigma})\\
&+(\nabla_h(\boldsymbol\sigma-\boldsymbol\sigma_h),
\nabla_h(\hat{\boldsymbol\sigma}-I_h^{\mathrm{VE}}\hat{\boldsymbol\sigma}))\\
&+(\operatorname{curl}_h(\boldsymbol\sigma-\boldsymbol\sigma_h),\hat{\boldsymbol p})\\
&-(\nabla\hat{\boldsymbol\sigma}\,\boldsymbol n+\hat{\boldsymbol p}\times\boldsymbol n,
\boldsymbol\sigma-\boldsymbol\sigma_h)_{\partial\mathcal T_h}.
\end{aligned}
\end{equation}
	
We estimate these terms separately. For \(k=0\), \(I_h^{\mathrm{VE}}\) is the CR interpolant and \(\widetilde{\hat{\boldsymbol\sigma}}=I_h^{\mathrm{VE}}\hat{\boldsymbol\sigma}\). Since the CR interpolant preserves the elementwise mean gradient,
\[
\begin{aligned}
(\nabla_h(\boldsymbol\sigma-\boldsymbol\sigma_h),
\nabla_h I_h^{\mathrm{VE}}\hat{\boldsymbol\sigma})
=
(\operatorname{grad}_h^\Pi(I_h^{\mathrm{VE}}\boldsymbol\sigma-\boldsymbol\sigma_h),
\operatorname{grad}_h^\Pi I_h^{\mathrm{VE}}\hat{\boldsymbol\sigma}).
\end{aligned}
\]
Applying the identity~\eqref{eq:dual_testing_identity}, the polynomial projection approximation estimates, and the CR interpolation estimate~\eqref{eq:VE_local_interpolation}, we obtain
\begin{equation}\label{tmp_k0_neg}
\begin{aligned}
~&\left|(\nabla_h(\boldsymbol\sigma-\boldsymbol\sigma_h),
\nabla_h I_h^{\mathrm{VE}}\hat{\boldsymbol\sigma})\right|\\
\lesssim~&
h^2(|\boldsymbol\sigma|_2+\|\boldsymbol g\|)|\hat{\boldsymbol\sigma}|_2
+h^2\|\boldsymbol g-\boldsymbol g_h\|\,|\hat{\boldsymbol\sigma}|_2
+\|\boldsymbol g-\boldsymbol g_h\|_{-1}\|\hat{\boldsymbol\sigma}\|_1.
\end{aligned}
\end{equation}
In addition, we have
\[
\begin{aligned}
(\operatorname{curl}_h(\boldsymbol\sigma-\boldsymbol\sigma_h),\hat{\boldsymbol p})
=
(\operatorname{curl}_h(\boldsymbol\sigma-\boldsymbol\sigma_h),
\hat{\boldsymbol p}-I_h^{\mathrm{nt}}\hat{\boldsymbol p})
+(\operatorname{dev}\operatorname{grad}_w\boldsymbol r_h,
I_h^{\mathrm{nt}}\hat{\boldsymbol p}).
\end{aligned}
\]
The interpolation estimate~\eqref{interpolation_error_p} with \(m=1\) and the argument used to prove \eqref{eq_err5} imply
\[
\left|(\operatorname{curl}_h(\boldsymbol\sigma-\boldsymbol\sigma_h),\hat{\boldsymbol p})\right|
\lesssim
h\|\hat{\boldsymbol p}\|_1
\bigl(|\boldsymbol\sigma-\boldsymbol\sigma_h|_{1,h}
+\interl{\boldsymbol r_h}\bigr).
\]

The second term in~\eqref{eq:prop_tmp} is bounded by the CR interpolation estimate \eqref{eq:VE_local_interpolation}:
\[
\begin{aligned}
\left|(\nabla_h(\boldsymbol\sigma-\boldsymbol\sigma_h),
\nabla_h(\hat{\boldsymbol\sigma}-I_h^{\mathrm{VE}}\hat{\boldsymbol\sigma}))\right|
&\le
|\boldsymbol\sigma-\boldsymbol\sigma_h|_{1,h}
|\hat{\boldsymbol\sigma}-I_h^{\mathrm{VE}}\hat{\boldsymbol\sigma}|_{1,h}\\
&\lesssim
h|\boldsymbol\sigma-\boldsymbol\sigma_h|_{1,h}
|\hat{\boldsymbol\sigma}|_2.
\end{aligned}
\]

For the boundary term, the vanishing face means of the jumps, the trace inequality, and the local Poincaré inequality yield
\[
\begin{aligned}
\left|(\nabla\hat{\boldsymbol\sigma}\,\boldsymbol n+\hat{\boldsymbol p}\times\boldsymbol n,
\boldsymbol\sigma-\boldsymbol\sigma_h)_{\partial\mathcal T_h}\right|
\lesssim 
h|\boldsymbol\sigma-\boldsymbol\sigma_h|_{1,h}
(\|\hat{\boldsymbol\sigma}\|_2+\|\hat{\boldsymbol p}\|_1).
\end{aligned}
\]
Combining these bounds with the energy error estimate~\eqref{eq:stokes_full_gradient}, the dual regularity estimate~\eqref{eq:k0_dual}, and \(\|\boldsymbol g\|\lesssim|\boldsymbol\sigma|_2+|\boldsymbol p|_1\), we obtain \eqref{eq:k0_L2}.
\end{proof}

\section{Application to a decoupled triharmonic scheme}\label{sec:triharmonic}

We apply the distributional Stokes method of Section~\ref{sec:discrete_generalized_stokes} to a decoupled scheme for the triharmonic equation~\eqref{triharmonic}, using \(H^2\)-nonconforming virtual elements for the two biharmonic subproblems. Suitable load and tensor reconstructions yield optimal Hessian error estimates under the stated regularity assumptions.


\subsection{Decoupled scheme}\label{sec:tensorreconstruct}
Let \(f\in L^2(\Omega)\) be the source term in the triharmonic problem~\eqref{triharmonic}. Following \cite{ChenHuang2018,GuHuang2026}, we decouple the triharmonic problem~\eqref{triharmonic} into two biharmonic problems and a curl–dev grad Stokes problem. Specifically, we seek \(w,u\in H_0^2(\Omega)\) and a solution \((\boldsymbol\sigma,\boldsymbol p,\boldsymbol r)\) of the continuous Stokes formulation with \(\boldsymbol g=\nabla^2w\), such that
\begin{subequations}\label{triharmonic_vem_decoupled}
\begin{align}
(\nabla^2w,\nabla^2v)&=(f,v), &&\forall\,v\in H_0^2(\Omega),\label{decoupled_triharmonic_continue1}\\
(\nabla^2u,\nabla^2\chi)&=(\bs\sigma,\nabla^2\chi), &&\forall\,\chi\in H_0^2(\Omega).\label{decoupled_triharmonic_continue4}
\end{align}
\end{subequations}
By Theorem~\ref{thm:gStokes_continue_wellposed}, this formulation is equivalent to the splitting in \cite{GuHuang2026}. In particular, the resulting \(u\) belongs to \(H_0^3(\Omega)\) and is the unique solution of the triharmonic problem~\eqref{triharmonic}, with \(\boldsymbol\sigma=\nabla^2u\).

For $ k \geqslant 0$, let $\Whvem$ be the $H^2$-nonconforming virtual element space of \cite{ChenHuang2020,ChenHuang2025a}, with interior degree $r=k-2$ in \cite{ChenHuang2025a} and homogeneous clamped boundary conditions, and let $\tIvem$ be its canonical interpolant. For $k=0$, this is the quadratic Morley--Wang--Xu space \cite{morley2006}. We recall the following norm equivalence \cite[(4.14)]{ChenHuang2025a}, commuting identity \cite[(4.3), (4.12)]{ChenHuang2025a}, and interpolation estimate \cite[Lemma~5.3]{ChenHuang2020}:
\begin{subequations}
\begin{align}
\|Q_{k,h}\nabla_h^2v_h\| &\simeq\|\nabla_h^2v_h\|, && v_h \in \Whvem, \label{norm_equivalence_H2} \\
Q_{k,h}\nabla_h^2(\tIvem v)&=Q_{k,h}\nabla^2v,
\quad &&v\in H_0^2(\Omega),  \label{h2_vem_commuting} \\
| v - \tIvem v |_{2,h} &\lesssim h^{m-2} | v |_{m,h}, \quad &&v\in H_0^2(\Omega) \cap H^m(\mathcal{T}_h), 2 \leqslant m \leqslant k + 3.  \label{h2_vem_interpolation_estimate}
\end{align}
\end{subequations}

We use the biharmonic discretization of \cite[(4.13)]{ChenHuang2025a} with a modified load to preserve the optimal accuracy of the decoupled scheme. Let \(\PT\) be the local \(H^2\) projector defined in \cite[Section~3.3]{ChenHuang2020}. Define
\[
R_{k+2,T}v
:=Q_{k-2,T}v+(I-Q_{k-2,T})\PT v,
\qquad v\in H^2(T),
\]
which preserves polynomial moments of \(v\) up to degree \(k-2\) on \(T\). Let \(R_{k+2,h}\) and \(\Ph\) denote the operators obtained by applying \(R_{k+2,T}\) and \(\PT\) on each element, respectively. We then replace the original load term by \((f,R_{k+2,h}v_h)\). For \(k=0\), \(R_{2,h}v_h=v_h\) on \(\mathring W_{2,h}^{\mathrm{VE}}\), so the load reduces to \((f,v_h)\).
\begin{lemma}\label{lem:Rh}
For $k\ge0$, the operator $ R_{k+2,h} $ satisfies
\begin{subequations}
\begin{align}
R_{k+2,h}\tIvem v &= R_{k+2,h}v, \qquad &&\forall\,v\in H_0^2(\Omega), \label{Rkh_interpolation_property}\\
(v - R_{k+2,h}v, q) &= 0, \qquad &&\forall\,v \in H^2(\mathcal{T}_h), q \in \mathbb{P}_{k-2}(\mathcal{T}_h), \label{Rkh_property1} \\
R_{k+2,h} v &= v , \quad && \forall\,v \in  \mathbb{P}_{k+2}(\mathcal{T}_h), \label{Rkh_property2}\\
\| v - R_{k+2,h} v \| &\lesssim h^{j} |v |_{j,h}, \qquad && \forall\,v \in H^{j}(\mathcal{T}_h), j=2,3.  \label{Rkh_property3}
\end{align}
\end{subequations}
\end{lemma}
\begin{proof}
Properties \eqref{Rkh_interpolation_property}--\eqref{Rkh_property2} follow from the defining moments, the projection--interpolation identity, and polynomial reproduction; see \cite[Lemma~3.7 and Corollary~3.10]{ChenHuang2020}. The \(H^2\)-seminorm stability of \(\PT\) \cite[(3.12)]{ChenHuang2020}, the Poincaré inequality on the kernel of \(\PT\) \cite[Lemma~A.5]{ChenHuang2020}, and the polynomial-preserving property of \(\PT\) \cite[Lemma~3.7]{ChenHuang2020}, together with standard polynomial approximation, yield
\[
\begin{aligned}
\|v-R_{k+2,T}v\|_T
&=\|(I-Q_{k-2,T})(v-\PT v)\|_T\le\|v-\PT v\|_T  \\
&\lesssim h_T^2|v-\PT v|_{2,T} \lesssim h_T^j|v|_{j,T}, v\in H^j(T),\quad j=2,3.
\end{aligned}
\]
Summing over the mesh proves \eqref{Rkh_property3}.
\end{proof}

For the tensor data in the final biharmonic block, we use the following standard reconstruction \cite[Section~3.4]{AyusoLipnikovManzini2016}.
For $T\in\mathcal T_h$, let $\Pi_{k+1,T}\colon H^1(T;\mathbb S)\to\mathbb P_{k+1}(T;\mathbb S)$ be the elliptic projector defined by
\begin{subequations}
\begin{align}
(\nabla\Pi_{k+1,T}\bs\sigma,\nabla\bs\tau)_T
&=(\nabla\bs\sigma,\nabla\bs\tau)_T,
\quad\forall\,\bs\tau\in\mathbb P_{k+1}(T;\mathbb S),\label{local_elliptic_1}\\
\int_T\Pi_{k+1,T}\bs\sigma\,dx&=\int_T\bs\sigma\,dx.\label{local_elliptic_2}
\end{align}
\end{subequations}
For discrete virtual tensors, this projector is computable from the DoFs \eqref{dof_vem}. Define $(\Pi_h\bs\sigma)|_T:=\Pi_{k+1,T}(\bs\sigma|_T)$. This operator reproduces $\mathbb P_{k+1}(\mathcal T_h;\mathbb S)$. For $k=0$, the discrete tensor is piecewise linear, hence $\Pi_h\bs\sigma_h=\bs\sigma_h$.
Standard stability and approximation properties give
\begin{subequations}
\begin{equation}\label{eq_rec_estimate1}
\|\Pi_h\bs\sigma\|\lesssim\|\bs\sigma\|+h|\bs\sigma|_{1,h},
\qquad\bs\sigma\in H^1(\mathcal T_h;\mathbb S).
\end{equation}
For $\bs\sigma\in H^\ell(T;\mathbb S)$ and $1\le\ell\le k+2$,
\begin{equation}\label{eq_rec_estimate3}
\|\bs\sigma-\Pi_{k+1,T}\bs\sigma\|_T
+h_T|\bs\sigma-\Pi_{k+1,T}\bs\sigma|_{1,T}
\lesssim h_T^\ell|\bs\sigma|_{\ell,T}.
\end{equation}
\end{subequations}

Combining biharmonic discretization with the Stokes method \eqref{discrete_generalized_stokes} gives the decoupled scheme: find
$(w_h,\bs\sigma_h,\bs p_h,\bs r_h,u_h) \in \Whvem \times \VhvemS \times \Vhnt \times(\VhRT/\mathbb R^3) \times \Whvemm$ such that
\begin{subequations}\label{decoupled_triharmonic_discrete}
\begin{align}
(Q_{k,h}\nabla_h^2w_h,Q_{k,h}\nabla_h^2v)& = (f, R_{k+2,h} v ),
&&\label{decoupled_triharmonic_discrete1}\\
a_h(\bs\sigma_h,\bs r_h;\bs\tau,\bs s)+b_h(\bs\tau,\bs s;\bs p_h)& = (Q_{k,h}\nabla_h^2w_h, \tilde{\bs\tau}),
&&\label{decoupled_triharmonic_discrete2}\\
b_h(\bs\sigma_h,\bs r_h;\bs q)&=0,
&&\label{decoupled_triharmonic_discrete3}\\
(Q_{k+1,h}\nabla_h^2u_h,Q_{k+1,h}\nabla_h^2\chi)&=(\Pi_h \bs\sigma_h,Q_{k+1,h}\nabla_h^2\chi),
&&\label{decoupled_triharmonic_discrete4}
\end{align}
for all $(v,\bs\tau,\bs q,\bs s,\chi)\in\Whvem\times\VhvemS\times\Vhnt\times(\VhRT/\mathbb R^3)\times \Whvemm $.
\end{subequations}
The norm equivalences on the two $H^2$-nonconforming virtual element spaces, together with Theorem~\ref{thm:gStokes_discrete_wellposed}, imply well-posedness of all three blocks.

\subsection{Error estimates}
Let $\delta_{k0}$ denote the Kronecker delta.

\begin{lemma}\label{lem:err_biharmonic_H2}
Let $k\ge0$, and let $w\in H_0^2(\Omega)\cap H^{k+3}(\Omega)$ and $w_h\in\Whvem$ solve \eqref{decoupled_triharmonic_continue1} and \eqref{decoupled_triharmonic_discrete1}, respectively. Then
\begin{equation} \label{err_biharmonic_H2}
\|Q_{k,h}\nabla_h^2(w-w_h)\|\lesssim h^{k+1}|w|_{k+3}+\delta_{k0}h^2\|f\|. 
\end{equation}
\end{lemma}
\begin{proof}
We follow the proof of \cite[Theorem~4.7]{ChenHuang2025a}, with the modified load-consistency term estimated as follows. For $k\ge1$ and $v_h\in\Whvem$, \eqref{Rkh_property1}, \eqref{Rkh_property3} with $j=2$, and $f=\Delta^2w$ give
\[
\begin{aligned}
|(f,v_h-R_{k+2,h}v_h)|
&=|(f-Q_{k-2,h}f,v_h-R_{k+2,h}v_h)|\\
&\lesssim h^{k+1}|w|_{k+3}|v_h|_{2,h}.
\end{aligned}
\]
For $k=0$, this term vanishes since $R_{2,h}v_h=v_h$. The remaining approximation and consistency estimates are unchanged, and \eqref{norm_equivalence_H2}--\eqref{h2_vem_commuting} yield \eqref{err_biharmonic_H2}.
\end{proof}

\begin{corollary}\label{thm_tri_stokes_error}
Let $k\ge0$. Let \((w,\boldsymbol\sigma,\boldsymbol p,0)\) satisfy the continuous biharmonic equation~\eqref{decoupled_triharmonic_continue1} and the continuous Stokes formulation~\eqref{continue_generalized_stokes} with \(\boldsymbol g=\nabla^2w\), and let \((w_h,\boldsymbol\sigma_h,\boldsymbol p_h,\boldsymbol r_h)\) solve the discrete equations \eqref{decoupled_triharmonic_discrete1}--\eqref{decoupled_triharmonic_discrete3}. Assume $w\in H^{k+3}(\Omega)$, $\bs\sigma\in H^{k+2}(\Omega;\bbS)$, and $\bs p\in H^{k+1}(\Omega;\bbT)$. Then
\begin{equation}\label{err_tri_stokes_H1}
\begin{aligned}
&\|\gradvem(\Ivem\bs\sigma-\bs\sigma_h)\|+\|\bs p-\bs p_h\|+\interl{\bs r_h}\\
&\qquad\lesssim h^{k+1}(|\bs\sigma|_{k+2}+|\bs p|_{k+1}+|w|_{k+3})+\delta_{k0}h^2\|f\|.
\end{aligned}
\end{equation}
Assume additionally that the dual regularity estimates~\eqref{dual_regularity} hold for \(k\ge1\), or that the dual estimate~\eqref{eq:k0_dual} holds for \(k=0\). Then
\begin{equation}\label{err_tri_stokes_L2}
\|\boldsymbol\sigma-\boldsymbol\sigma_h\|
\lesssim h^{k+2}
\bigl(|\boldsymbol\sigma|_{k+2}+|\boldsymbol p|_{k+1}+|w|_{k+3}\bigr)
+\delta_{k0}h^3\|f\|.
\end{equation}
\end{corollary}
\begin{proof}
For $\bs g=\nabla^2w$ and $\bs g_h=Q_{k,h}\nabla_h^2w_h$, the projection estimate and the biharmonic error estimate~\eqref{err_biharmonic_H2} give
\begin{equation}\label{corollary_eq1}
\|\bs g-\bs g_h\|\le\|(I-Q_{k,h})\nabla^2w\|+\|Q_{k,h}\nabla_h^2(w-w_h)\|
\lesssim h^{k+1}|w|_{k+3}+\delta_{k0}h^2\|f\|.
\end{equation}
Substituting this bound into the Stokes energy error estimate of Theorem~\ref{thm_H1bound} and using the gradient commuting identity~\eqref{com_grad} yield \eqref{err_tri_stokes_H1}.

For \(k\ge1\), we estimate \((\boldsymbol g-\boldsymbol g_h,\hat{\boldsymbol\sigma})\) in \eqref{eq_L2_tmp}  directly to avoid introducing a separate bound for \(\|\boldsymbol g-\boldsymbol g_h\|_{-1}\). Since \(\hat{\boldsymbol\sigma}=\nabla^2\hat u\), test the continuous biharmonic equation \eqref{decoupled_triharmonic_continue1} with \(\hat u\) and the discrete equation~\eqref{decoupled_triharmonic_discrete1} with \(\tIvem\hat u\). Using Lemma~\ref{lem:Rh}, the commuting identity~\eqref{h2_vem_commuting}, and the dual regularity estimate~\eqref{dual_regu}, we obtain
\[
\begin{aligned}
|(\boldsymbol g-\boldsymbol g_h,\hat{\boldsymbol\sigma})|
&=\left|(\nabla^2w,\nabla^2\hat u)
-(Q_{k,h}\nabla_h^2w_h,Q_{k,h}\nabla^2\hat u)\right|\\
&=|(f,\hat u-R_{k+2,h}\hat u)|
=|(f-Q_{k-2,h}f,\hat u-R_{k+2,h}\hat u)|\\
&\lesssim h^{k+2}|f|_{k-1}|\hat u|_3 \lesssim h^{k+2}|w|_{k+3}\|\boldsymbol e_T\|.
\end{aligned}
\]
Using this bound for \((\boldsymbol g-\boldsymbol g_h,\hat{\boldsymbol\sigma})\) in the decomposition~\eqref{eq_L2_tmp}, and the duality argument in Theorem~\ref{thm:stokes_L2_error}, we have 
\[
\| \bs \sigma - \bs \sigma_{h} \|
\lesssim h^{k+2}
\bigl(|\boldsymbol\sigma|_{k+2}+|\boldsymbol p|_{k+1}+|w|_{k+3}\bigr)
+h\|\boldsymbol g-\boldsymbol g_h\|.
\]
Combining this estimate with \eqref{corollary_eq1} then proves \eqref{err_tri_stokes_L2} for \(k\ge1\).

For \(k=0\), let \((\hat{\boldsymbol\sigma},\hat{\boldsymbol p})\) satisfy the dual problem~\eqref{eq:k0_dual}. Since \(\operatorname{curl}\hat{\boldsymbol\sigma}=0\), exactness of the Hessian complex~\eqref{hessiancomplex} implies there exists \(\hat u\in H_0^3(\Omega)\) such that
\[
\nabla^2\hat u=\hat{\boldsymbol\sigma},
\qquad
|\hat u|_3=|\hat{\boldsymbol\sigma}|_1
\lesssim\|\boldsymbol\sigma-\boldsymbol\sigma_h\|.
\]
Testing \eqref{decoupled_triharmonic_continue1} with \(\hat u\) and \eqref{decoupled_triharmonic_discrete1} with \(\tIvem\hat u\), and using \eqref{h2_vem_commuting}, \eqref{Rkh_interpolation_property}, and \eqref{Rkh_property3} with \(j=3\), we obtain
\[
\begin{aligned}
	|(\boldsymbol g-\boldsymbol g_h,\hat{\boldsymbol\sigma})|
	&=|(f,\hat u-R_{2,h}\hat u)| \lesssim h^3\|f\|\,|\hat u|_3
	\lesssim h^3\|f\|\,\|\boldsymbol\sigma-\boldsymbol\sigma_h\|.
\end{aligned}
\]
Replacing the estimate of the dual data pairing in \eqref{tmp_k0_neg} by the preceding bound in the proof of Proposition~\ref{prop:k0_errors} and using the error estimate~\eqref{corollary_eq1} yields
\[
\|\boldsymbol\sigma-\boldsymbol\sigma_h\|
\lesssim h^2(|\boldsymbol\sigma|_2+|\boldsymbol p|_1+|w|_3)
+h^3\|f\|,
\]
which proves \eqref{err_tri_stokes_L2} for \(k=0\).
\end{proof}

\begin{remark}
Using the load $(f,\Ph v)$ in \eqref{decoupled_triharmonic_discrete1} does not provide sufficient load consistency for the preceding proof of \eqref{err_tri_stokes_L2} when $k\ge2$.
\end{remark}

\begin{theorem}\label{thm:tri_u_error}
Under the assumptions of Corollary~\ref{thm_tri_stokes_error} and the corresponding dual regularity assumption for the chosen $k\ge0$, let $u\in H_0^2(\Omega)$ and $u_h \in \Whvemm $ satisfy \eqref{decoupled_triharmonic_continue4} and \eqref{decoupled_triharmonic_discrete4}, respectively. Then
\begin{equation}\label{err_triharmonic_u_H2}
\begin{aligned}
&\|Q_{k+1,h}\nabla_h^2(u-u_h)\|+\|\nabla^2u-\nabla_h^2u_h\|\\
&\qquad\lesssim h^{k+2}(|\bs\sigma|_{k+2}+|\bs p|_{k+1}+|w|_{k+3})+\delta_{k0}h^3\|f\|.
\end{aligned}
\end{equation}
\end{theorem}
\begin{proof}
Let $z_h$ be the canonical interpolant of $u$ in $\Whvemm$. Applying the commuting identity~\eqref{h2_vem_commuting} with $k$ replaced by $k+1$ gives $Q_{k+1,h}\nabla_h^2 z_h = Q_{k+1,h}\bs\sigma$. Using this identity and testing the discrete biharmonic equation~\eqref{decoupled_triharmonic_discrete4} with \(z_h-u_h\), we obtain
\begin{equation}\label{biharmonic_uh3}
\|Q_{k+1,h}\nabla_h^2(u-u_h)\|=\|Q_{k+1,h}\nabla_h^2(z_h-u_h)\|
\le\| \bs \sigma- \Pi_h \bs\sigma_h\|.
\end{equation}
The reconstruction stability estimate~\eqref{eq_rec_estimate1} and approximation estimate~\eqref{eq_rec_estimate3} give
\[
\begin{aligned}
\|\boldsymbol\sigma-\Pi_h\boldsymbol\sigma_h\|
&\le\|\boldsymbol\sigma-\Pi_h\boldsymbol\sigma\|
+\|\Pi_h(\boldsymbol\sigma-\boldsymbol\sigma_h)\|\\
&\lesssim h^{k+2}|\boldsymbol\sigma|_{k+2}
+\|\boldsymbol\sigma-\boldsymbol\sigma_h\|
+h|\boldsymbol\sigma-\boldsymbol\sigma_h|_{1,h}\\
&\lesssim h^{k+2}
\bigl(|\boldsymbol\sigma|_{k+2}+|\boldsymbol p|_{k+1}+|w|_{k+3}\bigr)
+\delta_{k0}h^3\|f\|,
\end{aligned}
\]
where the last inequality follows from \eqref{norm_sigma}, \eqref{eq:VE_local_interpolation}, and \eqref{err_tri_stokes_H1}--\eqref{err_tri_stokes_L2}.
Since \(\boldsymbol\sigma=\nabla^2u\), we have \(u\in H^{k+4}(\Omega)\) and \(|u|_{k+4}\lesssim|\boldsymbol\sigma|_{k+2}\). The interpolation estimate~\eqref{h2_vem_interpolation_estimate} and the norm equivalence \eqref{norm_equivalence_H2} on \(\Whvemm\) yield
\[
\begin{aligned}
\|\nabla^2u-\nabla_h^2u_h\|
&\le |u-z_h|_{2,h}+|z_h-u_h|_{2,h}\\
&\lesssim h^{k+2}|u|_{k+4}
+\|Q_{k+1,h}\nabla_h^2(z_h-u_h)\|.
\end{aligned}
\]
Combining these bounds with \eqref{biharmonic_uh3} proves \eqref{err_triharmonic_u_H2}. 
\end{proof}

\section{Numerical Experiments}\label{sec:numerical_results}
We report numerical results for the distributional Stokes method \eqref{discrete_generalized_stokes} and the decoupled triharmonic scheme \eqref{decoupled_triharmonic_discrete}, focusing on convergence and the effect of the Stokes load reconstruction.

\begin{example}\label{example2}
\rm
We test the convergence of \eqref{discrete_generalized_stokes} and the effect of the load reconstruction on $\Omega=(0,\pi)^3$. Choose the exact solution 
\begin{equation*}
\boldsymbol{\sigma} = \nabla^2(\sin x \sin y \sin z)^3,  
\qquad \boldsymbol{r} = 0,
\end{equation*}
\begin{equation*}
\boldsymbol{p} = \curl
\begin{pmatrix}
0 & \sin x \sin y \sin z & 0  \\
\sin x \sin y \sin z & 0 & 0  \\
0 & 0 & 0  
\end{pmatrix}.
\end{equation*}
The load is chosen as $\boldsymbol g=-\Delta\boldsymbol\sigma+\sym\curl\boldsymbol p$. These fields satisfy the constraints and the homogeneous boundary conditions in \eqref{generalized_stokes}. 

Table~\ref{table_error_stokes_optimal} reports the projected-gradient errors for \(k=0,1,2\) with the prescribed load term \((\bs g,\tilde{\bs\tau})\). For the tensor $L^2$-type error, we report the full error \(\|\bs\sigma-\bs\sigma_h\|\) for \(k=0\) and the computable element-moment error \(\|\bs\sigma_T-\bs\sigma_{h,T}\|\) for \(k=1,2\). The projected-gradient errors and the $L^2$ errors in $\bs p$ approach order \(k+1\), while the reported tensor $L^2$-type errors approach order \(k+2\), consistently with \eqref{eq:stokes_full_gradient}, \eqref{compulable_L2_bound}, and Proposition~\ref{prop:k0_errors}. Table~\ref{table_error_stokes_k1} shows the results for \(k=1\) with the load \((\boldsymbol g,\boldsymbol\tau_T)\). The cell-projection load produces element-moment errors with rates tending toward two, compared with approximately three for the reconstructed load \((\boldsymbol g,\widetilde{\boldsymbol\tau})\) in Table~\ref{table_error_stokes_optimal}.
\begin{table}[htbp]
\centering
\caption{Errors of \eqref{discrete_generalized_stokes} with the prescribed load term $(\bs g,\tilde{\bs\tau})$. Here $e_{\sigma,\nabla}:=\|\grad\bs\sigma-\grad_h^\Pi\bs\sigma_h\|$, and $e_\sigma:=\|\bs\sigma-\bs\sigma_h\|$ for $k=0$ and $\|\bs\sigma_T-\bs\sigma_{h,T}\|$ for $k\ge1$.}
\label{table_error_stokes_optimal}
\begin{tabular}{cccccccc}
\toprule
$k$ & $h$ & $e_{\sigma,\nabla}$ & rate & $e_\sigma$ & rate & $\|\bs p-\bs p_h\|$ & rate \\
\midrule
\multirow{4}{*}{0}
& $2^{-2}\pi$ &14.273  &   --   & 2.2660  & --     & 5.0212 &  --    \\
& $2^{-3}\pi$ &7.5634  & 0.9162 & 0.65712 & 1.7859 & 2.8985 & 0.7927 \\
& $2^{-4}\pi$ &3.8313  & 0.9812 & 0.17535 & 1.9059 & 1.4520 & 0.9973 \\
\midrule
\multirow{3}{*}{1}
& $2^{-2}\pi$ &5.5065  &   --   & 0.11233   & --     & 2.2807  & --     \\
& $2^{-3}\pi$ &1.5928  &  1.7895& 0.010888  & 3.3669 & 0.58449 & 1.9642 \\
& $2^{-4}\pi$ &0.41802 &  1.9300& 0.0014222 & 2.9364 & 0.13657 & 2.0974 \\
\midrule
\multirow{3}{*}{2}
& $2^{-2}\pi$ &1.6792  &   --   & 0.26907   & --     & 1.2426  & --     \\
& $2^{-3}\pi$ &0.20511 &  3.0333& 0.026344  & 3.3524 & 0.20248 & 2.6176 \\
& $2^{-4}\pi$ &0.023006&  3.1563& 0.0018687 & 3.8174 & 0.02895 & 2.8061 \\
\bottomrule
\end{tabular}
\end{table}

\begin{table}[htbp]
\centering
\caption{Errors of \eqref{discrete_generalized_stokes} for $k=1$
with load $(\bs g,\bs\tau_T)$. The error symbols are defined in Table~\ref{table_error_stokes_optimal}.}
\label{table_error_stokes_k1}
\begin{tabular}{ccccccc}
\toprule
$h$ & $e_{\sigma,\nabla}$ & rate & $e_\sigma$ & rate
& $\|\bs p-\bs p_h\|$ & rate \\
\midrule
$2^{-2}\pi$ & 7.9298  &   --   & 1.3650  & --     & 3.6492  & --     \\
$2^{-3}\pi$ & 2.5402  & 1.6423 & 0.51774 & 1.3987 & 1.3599  & 1.4241 \\
$2^{-4}\pi$ & 0.69400 & 1.8719 & 0.14757 & 1.8108 & 0.37713 & 1.8503 \\
\bottomrule
\end{tabular}
\end{table}

\end{example}

\begin{example}\label{example3}
\rm
Let $\Omega=(0,\pi)^3$ and choose the exact solution and load
\[
u(x,y,z)=\sin^3x\sin^3y\sin^3z,\qquad f=-\Delta^3u.
\]
We solve the decoupled scheme \eqref{decoupled_triharmonic_discrete1}--\eqref{decoupled_triharmonic_discrete4} with $k=1$ and compare two choices of the load term in \eqref{decoupled_triharmonic_discrete2}, keeping all other parts of the scheme unchanged.

Table~\ref{table_error_u} reports the projected Hessian error $\|Q_{2,h}\nabla_h^2(u-u_h)\|$. It converges at approximately third order with the reconstructed Stokes load $(Q_{1,h}\nabla_h^2w_h,\tilde{\bs\tau})$, consistently with Theorem~\ref{thm:tri_u_error}. With the cell-projection load $(Q_{1,h}\nabla_h^2w_h,\bs\tau_T)$, the observed rates tend toward two. This comparison illustrates the effect of the Stokes load reconstruction on the accuracy of the decoupled scheme.

\begin{table}[htbp]
\setlength{\abovecaptionskip}{0pt}
\setlength{\belowcaptionskip}{0pt}
\caption{Errors $\|Q_{2,h}\nabla_h^2(u-u_h)\|$ for two choices
of the Stokes load term in
\eqref{decoupled_triharmonic_discrete2}, with $k=1$.}
\label{table_error_u}
\centering
\begin{tabular}{ccccc}
\hline
$h$
& $(Q_{1,h}\nabla_h^2w_h,\tilde{\bs\tau})$
& rate
& $(Q_{1,h}\nabla_h^2w_h,\bs\tau_T)$
& rate \\
\hline
$2^{-2}\pi$ & 0.31830   & --     & 1.5537  & --     \\
$2^{-3}\pi$ & 0.051805  & 2.6192 & 0.52163 & 1.5746 \\
$2^{-4}\pi$ & 0.0063521 & 3.0278 & 0.14735 & 1.8238 \\
\hline
\end{tabular}
\end{table}
\end{example}

\section{Conclusion}\label{sec:conclusion}
We have developed an arbitrary-order distributional method for the curl--dev grad Stokes system through a tensor--multiplier coupling with an $H(\div)$ vector multiplier. The Raviart--Thomas moments support commuting interpolation, while the divergence-augmented norm equivalence controls the full multiplier space modulo translations; together, they give mesh-uniform stability for each fixed order. Under the stated solution regularity and load approximation assumptions, the energy error is optimal, and optimal full tensor $L^2$ accuracy additionally requires the dual regularity assumptions. In the triharmonic application, moment-preserving loads and tensor reconstruction transfer this accuracy through the decoupled scheme to the final Hessian approximation.

\bibliographystyle{abbrv}
\bibliography{ref}

\end{document}